\documentclass[preprint,1p,times]{elsarticle}

\usepackage{amssymb,amsmath}
\usepackage{amsthm}
\usepackage{color}
\usepackage{ulem}
\usepackage{graphicx}
\usepackage{latexsym}
\usepackage{mathptmx}

\newcommand{\bm}[1]{\bf #1}
\newcommand{\eps}{\epsilon}
\newcommand{\pa}{\partial}
\newcommand{\rf}[1]{(\ref{#1})}

\newcommand{\ddfrac}[2]{\displaystyle \frac{#1}{#2}}
\newcommand{\dint}{\displaystyle \int}
\newcommand{\dsum}{\displaystyle \sum}

\newcommand{\dlim}{\displaystyle \lim}
\newcommand{\dmax}{\displaystyle \max}

\newcommand{\vep}{\varepsilon}

\newcommand{\ode}[2]{\displaystyle \frac{d #1}{d #2}}

\newcommand{\ovl}[1]{\overline{#1}}

\newcommand{\wt}[1]{\widetilde{#1}}
\newcommand{\wh}[1]{\widehat{#1}}
\newcommand{\Qed}
  {\hfill\hbox{\rule[-2pt]{3pt}{6pt}}\par\bigskip\noindent}

\newcommand{\SECT}[1]{\setcounter{equation}{0}\section{#1}\quad}

\newtheorem{PROP}{Proposition}[section]
\newtheorem{REM}{Remark}[section]

\begin{document}

\begin{frontmatter}



\title{
Effective nonlocal kernels on Reaction-diffusion networks
\tnoteref{t1}
}
\tnotetext[t1]{
This work was supported in part by JST CREST (No.~JPMJCR14D3) to S.-I.~E., 
JSPS KAKENHI (No. 17K14228) to Y. T.
}

\author[hokudai]{
Shin-Ichiro~Ei\corref{cor}}
\ead{Eichiro@math.sci.hokudai.ac.jp}

\author[hokudai]{
Hiroshi~Ishii}
\ead{hiroshi-ishii@eis.hokudai.ac.jp}

\author[osaka]{
Shigeru~Kondo}
\ead{skondo@fbs.osaka-u.ac.jp}

\author[kyushu]{
Takashi~Miura}
\ead{miura\_t@anat1.med.kyushu-u.ac.jp}

\author[hakodate]{
Yoshitaro~Tanaka}
\ead{yoshitaro.tanaka@gmail.com}

\address[hokudai]{ 
Department of Mathematics, 
Faculty of Science, 
Hokkaido University
}

\address[osaka]{Graduate School of Frontier Biosciences, Osaka University}

\address[kyushu]{
Graduate School of Medical Sciences, Kyushu University
}
\address[hakodate]{
School of Systems Information Science, Future University Hakodate
}

\cortext[cor]{Corresponding author. E-mail: Eichiro@math.sci.hokudai.ac.jp}

\begin{abstract}
A new method to derive an essential integral kernel from any given 
reaction-diffusion network is proposed. Any network describing
metabolites or signals with arbitrary many factors
can be reduced to a single or a simpler
system of integro-differential equations called ``effective equation''
including the reduced integral kernel (called ``effective kernel'' )
in the convolution type.
As one typical example, the Mexican hat shaped kernel is theoretically
derived from two component activator-inhibitor systems.
It is also shown that a three component system with quite different
appearance from activator-inhibitor systems is reduced to an
effective equation with the Mexican hat shaped kernel.
It means that the two different systems have essentially the same 
effective equations and
that they exhibit essentially the same spatial and temporal patterns.
Thus, we can identify two different systems with
the understanding in unified concept through the reduced effective kernels.
Other two applications of this method are also given:
Applications to pigment patterns on skins 
(two factors network with long range interaction) and
waves of differentiation (called proneural waves)
 in visual systems on brains (four factors network 
with long range interaction) . 
In  the applications, we observe the reproduction of
the same spatial and temporal patterns as those appearing in
pre-existing models through the numerical simulations of the effective equations.
\end{abstract}


\end{frontmatter}




\SECT{Introduction}

\begin{figure}  \label{fig1}
 \centering
  \includegraphics[width= 10cm, bb=0 0 612 230 ]{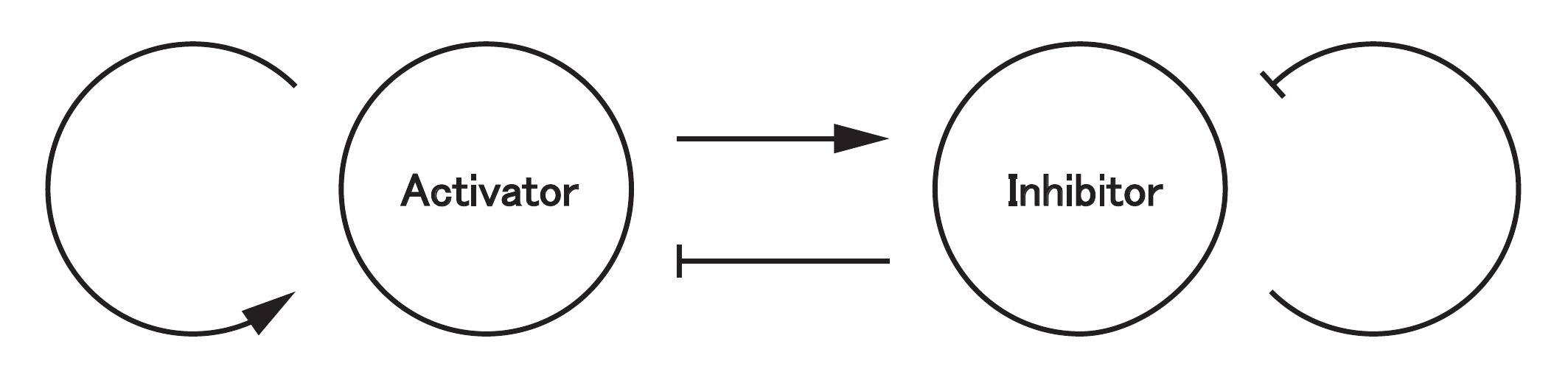}
 \caption{
 Activator-Inhibitor network. $\rightarrow$ and $\dashv$ denote activation
 and inhibition, respectively.
 }
\end{figure}

\begin{figure}[h]  
 \centering
  \includegraphics[width= 5cm, bb=0 0 429 320 ]{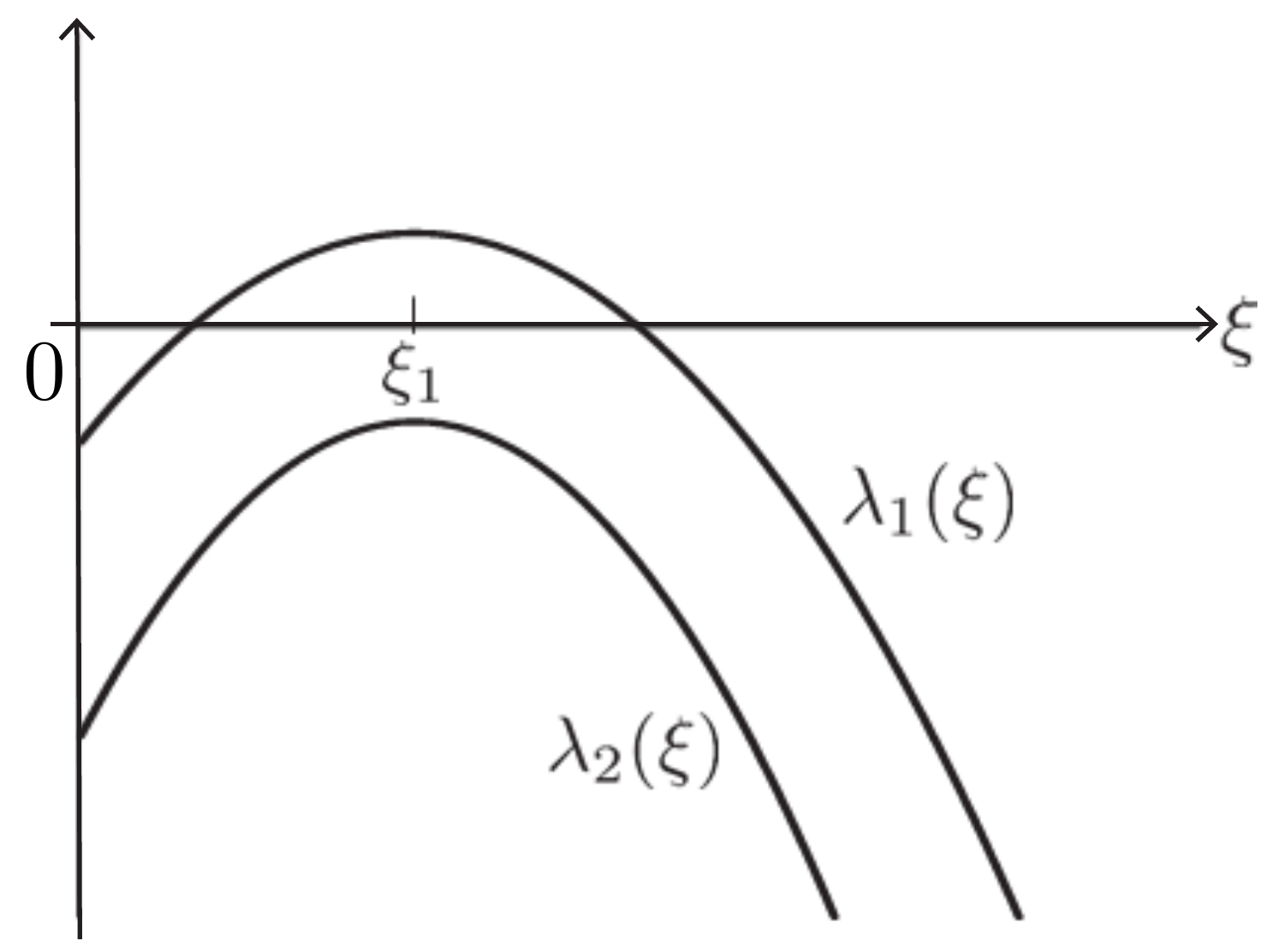}
 \caption{
 Schematic graphs of $\lambda_1(\xi )$ and $\lambda_2(\xi )$.
 In this case, $\lambda_{max}(\xi ) = \lambda_1(\xi )$ holds.
 }
 \label{fig2}
\end{figure}

The understanding of the mechanism for self-organization is one of
the  most attractive and important themes  in biology.
Diffusion driven instability proposed by Turing \cite{Tu} has given
a significant part of the  theoretical basement
which has been adapted by many works related to pattern formation problems
(Meinhardt \cite{Ma}, Murray \cite{Mu}).
In practice, the existence of the mechanism in real nature
has been observed and checked for several real phenomena (e.g. Ouyang and Swinney \cite{OS},
Kondo and Asai \cite{KA}, Yamaguchi et.al. \cite{YYK}).
The mechanism of the diffusion driven instability
was understood as the interaction between the slow diffusivity of short range 
self-enhancing activator and
 the fast one of long range inhibitor or simply
`` local activation with long-range inhibition (LALI)''
  (e.g. Gierer-Meinhardt \cite{GM},
 Oster \cite{Ost}) as in Fig.\ref{fig1}.
 The situation is frequently modeled in the form of reaction-diffusion systems with
 two components of the activator 
 $u=u(t,x)$ and the inhibitor $v=v(t,x)$ at  time $t>0$ and position $x$
 \begin{equation} \label{RD1}
 \left\{ \begin{array}{ccl}
 u_t &=& d_1\Delta u + f(u,v), \\
 v_t &=& d_2\Delta v + g(u,v),
 \end{array} \right.
 \end{equation}
 which is called Activator-Inhibitor system.
In \rf{RD1}, $\Delta$ denotes the Laplace operator and it is assumed that 
 the diffusion coefficients $0 < d_1 < d_2$ and the nonlinear terms
 $f_u > 0$, $f_v < 0$, $g_u > 0$ 
 and $g_v < 0$. The typical example of $f$ and $g$
 are $f(u,v) = c_1u - c_2v$ and $g(u,v) = c_3u - c_4v$ 
 with all positive constants $c_j, (j=1,\cdots, 4)$,
 which appear as a linear approximation in the neighbourhood of an equilibrium.
 In fact,  the diffusion driven instability is understood for the linearized system as follows:
 Let the linearized system for an equilibrium, say 0, be
 \begin{equation} \label{RD2}
 \left\{ \begin{array}{ccl}
 u_t &=& d_1\Delta u + c_1u - c_2v, \\
 v_t &=& d_2\Delta v + c_3u - c_4v.
 \end{array} \right.
 \end{equation}
 For the simplicity, we consider \rf{RD2} on $\bm{R}^1$. Then the Fourier
 transformation of \rf{RD2} yields
 \begin{equation} \label{RD3}
 \left\{ \begin{array}{ccl}
 \wh{u}_t &=& -d_1\xi^2\wh{u} + c_1\wh{u} - c_2\wh{v}, \\
 \wh{v}_t &=& -d_2\xi^2\wh{v} + c_3\wh{u} - c_4\wh{v}
 \end{array} \right.
 \end{equation}
while the specific notations are defined in Section 2.
 That is, \rf{RD3} is written as
 $\wh{U}_t = B(\xi )\wh{U}$ with 
 $\wh{U} := \left(\begin{array}{c} \wh{u} \\
 \wh{v}\end{array}\right
 )$
 and $B(\xi ) := \left( \begin{array}{cc} -d_1\xi^2 + c_1 & -c_2 \\
 c_3 & -d_2\xi^2 - c_4 \end{array}\right)$. Denoting the eigenvalues
 of $B(\xi )$ by $\lambda_1 (\xi )$, $\lambda_2(\xi )$ with
 $\lambda_1 (0) \ge \lambda_2(0)$ and defining $\lambda_{max}(\xi )
 := \max\{\lambda_1(\xi ),\lambda_2(\xi )\}$,
 we assume
 the maximal eigenvalue $\lambda_{max}(\xi )$ satisfies 
 $\lambda_{max}(0) < 0$ and $\lambda_{max}(\xi_1) > 0$ for
 $\xi_1 > 0$ as in Fig.\ref{fig2}. When the assumption for $\lambda_{max}(\xi )$ holds,
it is called ``the diffusion driven instability'' or ``Turing Instability''. 

\begin{figure}[h]  
 \centering
 \includegraphics[width=12cm,  bb=0 0 824 291]{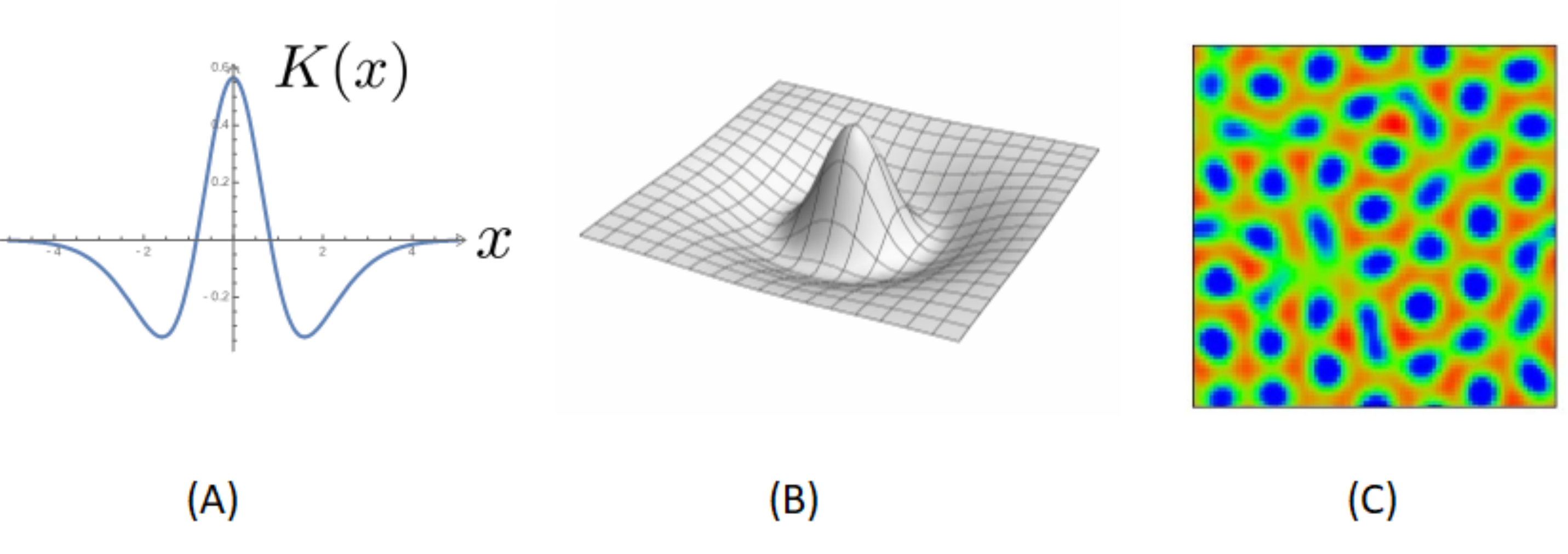}
 \caption{
 Integral kernels $K$ with the Mexican hat profile.
 (A) is on ${\bf R}$, (B) is on ${\bf R}^2$, 
 (C)  is the numerically
generated 2D patterns with Mexican hat kernel of (B). 
 }
 \label{fig3}
\end{figure}

 Recently, Kondo \cite{Kon1} suggested that the activator-inhibitor system
 is equivalent to a kernel-based model with the Mexican hat kernel  
 which is a function defined
 by regarding locally positive and  spreading negative parts as
 a local activation and a lateral inhibition, respectively as in Fig.\ref{fig3} 
 for $x \in \bm{R}$ or $\bm{x} \in \bm{R}^2$. 
Kondo proposed the following model with convolution
(for the definition, see Section 2)
by using an integral kernel $K$
 \begin{equation} \label{KT1}
 u_t = \chi (K*u) - \alpha u,
 \end{equation} 
 where $\chi (r) := \left\{ \begin{array}{cc}
 0, & r \le 0, \\
 r, & 0 < r \le r^*, \\
 r^*, & r \ge r^*
 \end{array} \right. $ for a positive constant $r^*$
 and the integral kernel function
 $K = K(|\bm{x}|)$ is a radially symmetric function
 with Mexican hat profiles or other general functions.
The Mexican hat shaped kernel $K$
was given in \cite{Kon1} as the sum of two Gaussian functions and numerically
showed the self-organization of
spatially periodic 2D patterns similar 
to those appearing in the activator-inhibitor reaction diffusion systems
can be reproduced in \rf{KT1}.
In practice, the Fourier transformation of the kernel leads the dispersion relation
of \rf{KT1} and the existence of a unstable nonzero wave number is
easily checked, that is, it has a similar structure of Turing Instability.
In that sense, the kernel-based  model like \rf{KT1} is called
``Kernel-based Turing model (KT model)'' in Kondo (\cite{Kon1}). 

\begin{figure}[ht]  
 \centering
 \includegraphics[width=8cm,  bb=0 0 525 335]{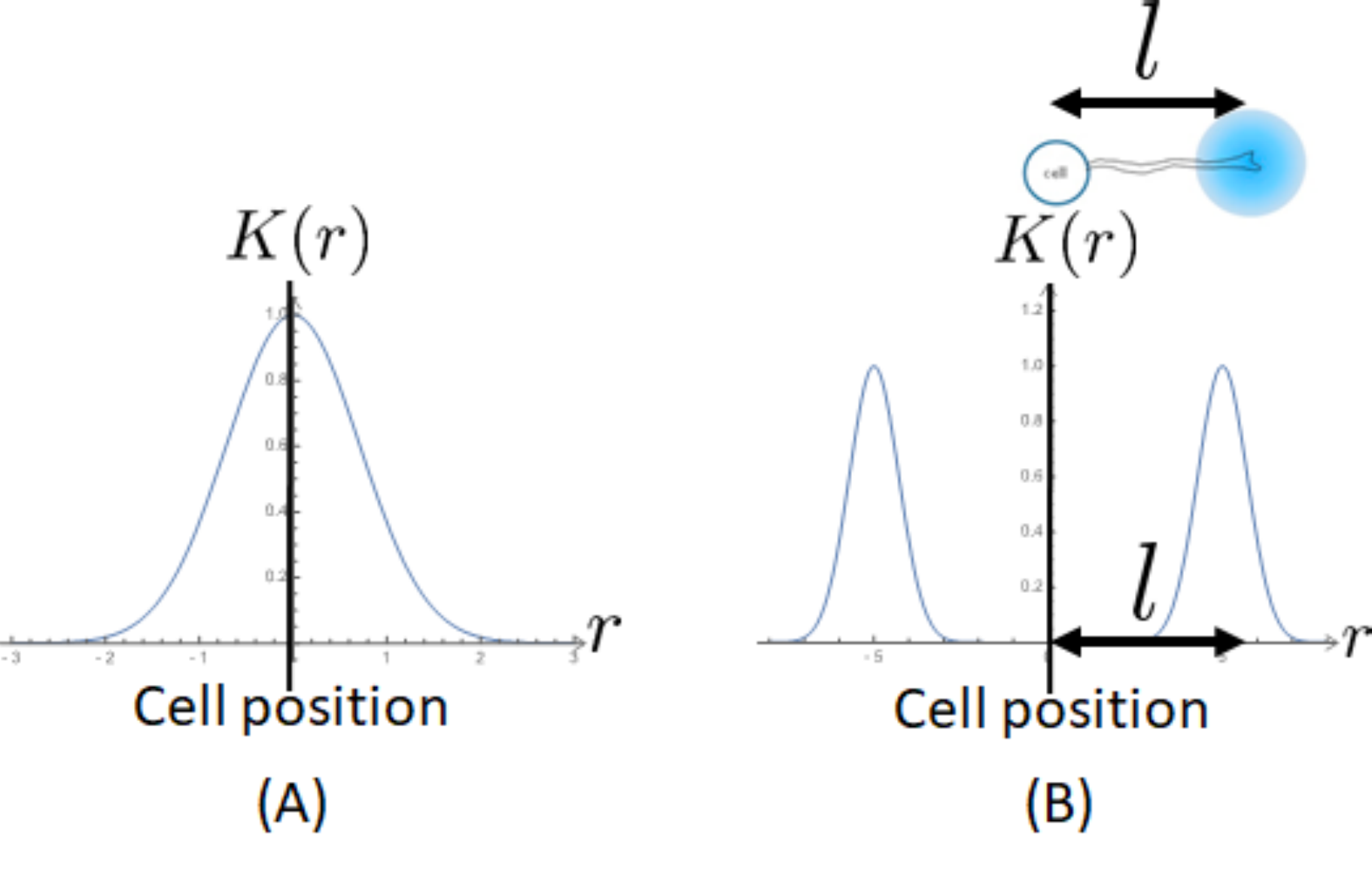}
 \caption{
 Examples of kernel shapes. (A) corresponds to diffusion process
 with the peak at the center ($r=0$). (B) is the kernel shape
by cell projections with peaks $l$ distant from the center.
 }
 \label{fig4}
\end{figure}

KT models in the type of \rf{KT1} are very effective to investigate 
spatially appearing patterns.
One reason is that we do not need to know
underlying mechanisms of molecules or cells in detail and that the
kernel shape is detected directly from an experimental observation
as stated in Kondo \cite{Kon1}.
 In fact, Kuffler \cite{Ku}  detected the kernel shape related to mammalian retina 
directly from the observation in real experiments.
For example, the diffusion process
is phenomenologically expressed by a unimodal 
localized profile of the kernel with one peak at the center as in Fig.\ref{fig4}(A),
whose relation was rigorously proved in Bates et. al. \cite{BC, BCC},
and the cell projection
which releases the signal molecule at the specific  position is
 by a kernel of a profile with peaks $l$ distant from the center
 as in Fig.\ref{fig4}(B).

Another reason is that
KT models can
easily reproduce more complicated patterns which are difficult to show by
conventional two component reaction diffusion systems such as nested patterns
appearing often on animal skins and sea shells. 
Such complicated patterns were numerically demonstrated in \cite{Kon1}
by using several other types of kernel functions.
Thus, spatial patterns and kernel shapes are directly related while
kernel shapes are related to the background mechanisms of
molecular or cellular processes 
together with the signal or metabolic networks. It means the possibility of 
two step approach
toward the theoretical understanding of spatial patterns, that is, one  is the
research into the
relation between patterns and kernel shapes which
corresponds to the macroscopic researches, another is between kernel shapes
and underlying biological mechanisms which corresponds to microscopic ones.  
Thus, KT models can connect macroscopic structures and microscopic ones
through the kernels and this two step approach
has a big advantage because these two steps can be treated independently
and separately.

Model equations of the type of \rf{KT1}  have been proposed
in many fields such as neural system (e.g. Amari\cite{Am}, Murray\cite{Mu} ),
cell-cell adhesion  problem (e.g. Carrillo et. al.\cite{CMSTT}, 
Painter et. al.\cite{PBSG} ),
Optical illusion (Sushida et. al. \cite{SKSM} ). 
Recently, integral kernels were also used to the development of 
the continuous method
for spatially discrete models (\cite{ESTY}).

In all of them,
appropriate integral kernels were adopted from the phenomenological point 
of view such as the direct detection from the experiments (e.g. Kuffler\cite{Ku}).
In that sense, they can be regarded as works in the step
for the relation between patterns and kernel shapes.
On the other hand, there have been almost no works in the step
between kernel shapes and underlying biological mechanisms
while both steps are absolutely necessary to understand the whole mechanism
from micro to macro  structures.

In this paper, we consider the relation between kernel shapes and 
the local networks of  metabolites or signals as underlying microscopic
mechanisms, 
and propose a new method to derive an effective kernel shape from
an arbitrarily given network.
Since the spatial and temporal patterns are essentially governed by
the effective equation with the effective kernel,
it can give us a unified aspect through the derived effective kernel
between seemingly different network systems.
Actually, we can identify seemingly different network systems
when the effective equations are same.

As one typical example of our approach, let us consider the activator-inhibitor system
with the local network like Fig.\ref{fig1}, which is basically described by the
reaction diffusion system \rf{RD1}.
As stated in Meinhardt and Gierer (\cite{GM}), the essential effect of this system is
the property of ``local activation with long-range inhibition (LALI)" and 
it has been believed that the corresponding kernel shape is the Mexican hat profile
as in Fig.\ref{fig3}. But there have been no theoretical investigation between
the kernel shape of the Mexican hat profile and the activator-inhibitor system \rf{RD1}. 
 By applying our technique, we systematically derive the Mexican hat shaped kernel 
 from \rf{RD2}. The details are shown in Section 3 together with the basic idea
 of this technique.
 
 We also demonstrate our technique  in Section 5 by applying to models
 with complicated networks
 including 3 components models in \cite{MDSM}, pigment patterns on skins
 with projections  in \cite{NTKK} and regulating waves 
 of differentiation in \cite{SYMMN}. 
 
 For the convenience of readers, we summarize notations and results
 of the Fourier transformation in Section 2.

 \SECT{Preliminaries of the Fourier transformation}

 In this section, 
 we give several notations and definitions for the Fourier transformation as follows:\\
 For a function $f(x)$ on $\bm{R}$, the Fourier transformation of $f$
 is defined by 
 \[\hat{f}(\xi ) = ({\cal F}f)(\xi )
 := \dint_{-\infty}^{\infty}f(x)e^{-ix\xi}dx\]
  and
 the inverse Fourier transformation is
 \[
 \check{f}(x) = ({\cal F}^{-1}f)(x) 
 := \ddfrac{1}{2\pi}\dint_{-\infty}^{\infty}f(\xi )e^{i\xi x}d\xi.
 \]
 In particular, 
 \begin{equation} \label{RD3-1}
 ({\cal F}^{-1}f)(x) 
 = \ddfrac{1}{\pi}\dint_0^{\infty}f(\xi^2 )\cos\xi xd\xi  
 \end{equation}
 holds for $f = f(\xi^2)$.
 The convolution of $f$ and $g$ is defined by
 \[
 (f*g)(x) := \dint_{-\infty}^{\infty}f(x-y)g(y)dy 
 = \dint_{-\infty}^{\infty}f(y)g(x-y)dy\]
  and
 $\wh{f*g} = \hat{f}\cdot\hat{g}$ holds. 

For a function $f(x,y)$ on $\bm{R}^2$, the Fourier transformation
  is defined by \[
  \hat{f}(\xi, \eta ) = ({\cal F}f)(\xi ,\eta )
 := \dint_{-\infty}^{\infty}\!\dint_{-\infty}^{\infty}
 f(x,y)e^{-i(x\xi+y\eta )}dxdy\]
  and
 the inverse Fourier transformation is
 \[\check{f}(x,y) = ({\cal F}^{-1}f)(x,y) 
 := \ddfrac{1}{(2\pi)^2}\dint_{-\infty}^{\infty}\!\dint_{-\infty}^{\infty}
 f(\xi ,\eta )e^{i(\xi x + \eta y)}d\xi d\eta.
 \]
 The convolution of $f$ and $g$ is
 \[
 (f*g)(x,y) := 
 \dint_{-\infty}^{\infty}\dint_{-\infty}^{\infty}f(x-x',y-y')g(x',y')dx'dy' 
 \] and
$\wh{f*g} = \hat{f}\cdot\hat{g}$ holds. 

 In particular, for a radially symmetric function
 $K = K(r)$ with $r := \sqrt{x^2+y^2}$, 
the Fourier transformation of $K$ is computed by polar coordinate
 \begin{eqnarray*}
 \hat{K}(\xi ,\eta ) 
 &=& \dint_{-\infty}^{\infty}\!\dint_{-\infty}^{\infty}
 K(\sqrt{x^2+y^2})e^{-i(\xi x + \eta y)}dxdy \\
 &=& \dint_0^{\infty}\!\dint_0^{2\pi}
 rK(r)e^{-ir(\xi\cos\theta + \eta\sin\theta )}drd\theta \\
 &=& \dint_0^{\infty}\!\dint_0^{2\pi}
 rK(r)e^{-irR\cos (\theta - \alpha )}drd\theta \\
 &=& \dint_0^{\infty}\!\dint_0^{2\pi}
 rK(r)e^{-irR\cos\theta }drd\theta,
 \end{eqnarray*}
 where $R := \sqrt{\xi^2+\eta^2}$ and $\alpha = \alpha(\xi ,\eta )$. That is,
 $\hat{K}$ is also radially symmetric and computed as
 \begin{eqnarray} 
 \nonumber
 \hat{K}(R) &=& 
\dint_0^{\infty}\!\!\!\!\dint_0^{2\pi}
 rK(r)e^{-irR\cos\theta }drd\theta \\
\nonumber
&=& 
2\dint_0^{\infty}\!\!\!\!\dint_0^{\pi}
 rK(r)\cos(rR\sin\theta )d\theta dr  \\
\label{e11-1} 
 &=&
4\!\dint_0^{\infty}\!\!\!\!\dint_0^{\pi/2}\!\!\!\!\!\!
 rK(r)\cos ( rR\sin\theta )d\theta dr .
\end{eqnarray}
 Similarly, the inverse Fourier transformation for a radially symmetric function $K(R)$ is
 \begin{eqnarray} 
 \nonumber
 \check{K}(r) &=& \ddfrac{1}{(2\pi)^2}
\dint_0^{\infty}\!\!\!\!\dint_0^{2\pi}
 RK(R)e^{irR\cos\theta }dRd\theta \\
\nonumber
&=& \ddfrac{1}{2\pi^2}
\dint_0^{\infty}\!\!\!\!\dint_0^{\pi}
 RK(R)\cos(rR\sin\theta )d\theta dR  \\
\label{e11-2} 
 &=&
 \ddfrac{1}{\pi^2}
 \dint_0^{\infty}\!\!\!\!\dint_0^{\pi/2}\!\!\!
 RK(R)\cos (rR\sin\theta )d\theta dR .
\end{eqnarray}

Here we give approximations of Dirac $\delta$ function, say
$\delta (x)$. In general, the following proposition holds.
\begin{PROP} (Approximation of $\delta (x)$) \label{prop1} \\
Suppose $G_{\eps}(x)$ is a function satisfying \\
1)~ $G_{\eps}(x) \ge 0$, \\
2)~ $\dint_{\bm{R}}G_{\eps}(x)dx = 1$,\\
3)~ $\dlim_{\eps \downarrow 0}\dint_{|x|>\eta}G_{\eps}(x)dx = 0$
for any $\eta > 0$.\\
Then $(G_{\eps}*u)(x) \rightarrow u(x)$ as $\eps \downarrow 0$ holds, 
that is, $G_{\eps} \rightarrow \delta (x)$.
\end{PROP} 
Typical examples of $G_{\eps}$ are the heat kernel  
$H_{\eps}(x) := \ddfrac{1}{\sqrt{4\pi\eps}}
e^{-\frac{x^2}{4\eps}}$, and a mollifier 
$J_{\eps}(\xi )$ given by $J_{\eps}(x ) := \frac{1}{\eps}J(\frac{x }{\eps})$
for a function $J(x)$ satisfying $J(x) \ge 0$, $\dint_{\bm{R}}J(x)dx = 1$ and
$J(x) = 0$ for $|x| \ge 1$.
In this paper, we furthermore 
impose the condition $\wh{J}(\xi ) \ge 0$
in order to keep the order of eigenvalues while the condition is not required
in general. Such a function $J(x)$ is for example made by
$J(x) = \{j*j\}(x)$ for a function $j(x)$. Then 
$\wh{J}(\xi ) = \{\wh{j}(\xi )\}^2 \ge 0$ holds
while $\wh{H_{\eps}}(\xi ) = e^{-\eps\xi^2} > 0$ holds.
Since the mollifier $J_{\eps}(x)$ satisfies
$(J_{\eps}*u)(x) = 0$ outside of $\eps-$ neighborhood of the support of
$u(x)$, that is,
$supp (J_{\eps}*u) \subset W_{\eps}(supp (u))$ holds
for $W_{\eps}(supp(u)) := \{x+y;\;
x \in supp(u),\; |y| \le \eps\} $,
it is useful when the support should be taken into account.
Above properties for $\delta (x)$
will be used to compute several kernels in the following sections.
 
 In the last of this section, we give several properties of the Lambert W function,
 say $w = W(z)$, which is a root $w$ of the equation
 \begin{equation} \label{Lam}
 we^w = z.
 \end{equation}
 Then we have:
 \begin{PROP}  \label{propLW}
 
 \ \\
 1)~ For $z \ge 0$, \rf{Lam} has a non-negative root $w_0(z) \ge 0$ and any other root
 $w$ satisfies $Re (w) < 0$.
 \\
 2)~ For $-1/e < z < 0$, \rf{Lam} has two real roots $w_1(z)$, $w_2(z)$ 
 ($0 >w_1(z) > -1 > w_2(z)$)
 and any other root $w$ satisfies $Re(w) < w_1(z)$. 
 \\
 3)~ For $z = -1/e$, \rf{Lam} has the root $w_0 = -1$ and any other root
$w$ satisfies $Re(w) < -1$.
 \end{PROP}
 This propositon is referred to \cite{CGHJK}.
 Here we express the principal branch of the Lambert W function
 by $W_0(z)$, which is given by 
 $W_0(z) = w_0 \in {\bf C}$ satisfying $Re (w_0) = \max\{ Re(w);\; we^w = z\}$.
 It is noted by Proposition \ref{propLW}
 that $W_0(z)$ is real for $z \ge -1/e$ satisfying
 $W_0(z) = w_0(z)$ for $z \ge 0$ in 1), $W_0(z) = w_1(z)$ 
 for $-1/e < z < 0$ in 2) and $W_0(-1/e) = -1$ in 3) of Proposition \ref{propLW}.
 
 \SECT{The basic idea
 and the theoretical derivation of the Mexican hat type kernel}

\begin{figure}[h]  
 \centering
   \includegraphics[width= 10cm, bb=0 0 609 195 ]{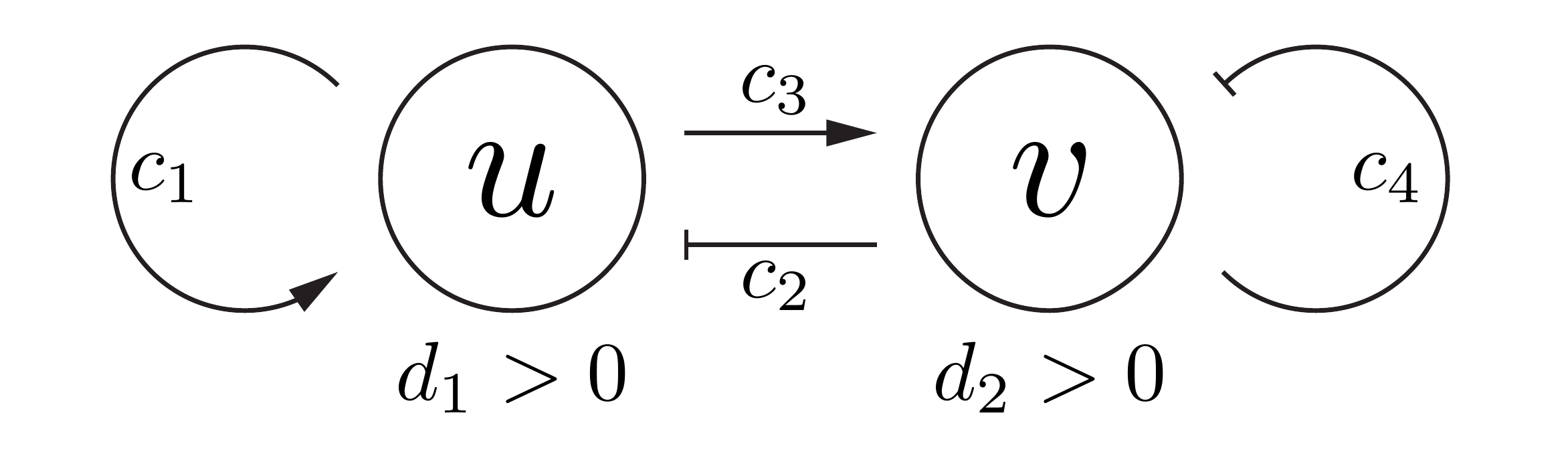}
 \caption{
 The network of Activator-Inhibitor system. $\rightarrow$ and
 $\dashv$ denote activation and inhibition, respectively.
 }
 \label{fig5}
\end{figure}

 In this section, we show our basic idea by using the activator-inhibitor 
 system \rf{RD2} (Fig.\ref{fig5})
 and theoretically derive the Mexican hat type kernel from it. For simplicity,
  we consider it on ${\bf R}$.
  
  The activator-inhibitor system has been represented by reaction diffusion systems,
  which have the following form in a linearized sense
  \begin{equation} \label{RD3-2}
 \left\{ \begin{array}{ccl}
 u_t &=& d_1u_{xx} + c_1u - c_2v, \\
 v_t &=& d_2v_{xx} + c_3u - c_4v
 \end{array} \right.
 \end{equation}
 for $x \in \bm{R}$, where
 all coefficients $d_1$, $d_2$ and $c_j (j=1,\cdots,4)$
 are positive and $d_1 < d_2$.
\rf{RD3-2} is written in a vector valued
 form 
 \begin{equation} \label{RD4}
 U_t = DU_{xx} + AU,
 \end{equation} 
 where $U = U(t,x) := \left(\begin{array}{c}
 u(t,x) \\ v(t,x)\end{array}\right) \in \bm{R}^2$,
 $D := diag\{d_1, d_2\}
 = \left( \begin{array}{cc} d_1 & 0 \\
 0 & d_2\end{array}\right)$ and
 $A := \left( \begin{array}{cc} c_1 & -c_2 \\
 c_3 & -c_4 \end{array} \right)$.
Here we note that variables $u$ and $v$ can take any sign
since we consider solutions in the neighborhood of an equilibrium by
the linearized equation like \rf{RD3-2}.

Taking the Fourier transformation of \rf{RD4}, we get
\begin{equation} \label{RD4-1}
\wh{U}_t = B(\xi )\wh{U} ,
\end{equation}
where
$
B(\xi ) := -\xi^2D + A
= \left(
\begin{array}{cc}
-d_1\xi^2 + c_1 & -c_2 \\
c_3 & -d_2\xi^2 - c_4
\end{array} \right) $.
Let $\lambda_j(\xi )$ ($j = 1, 2$)  be eigenvalues of $B(\xi )$ and
define $\lambda_{max}(\xi ) := \dmax_j\{ \lambda_j(\xi )\}$.
\begin{figure} 
 \centering
 \includegraphics[width=7cm,  bb=0 0 462 348]{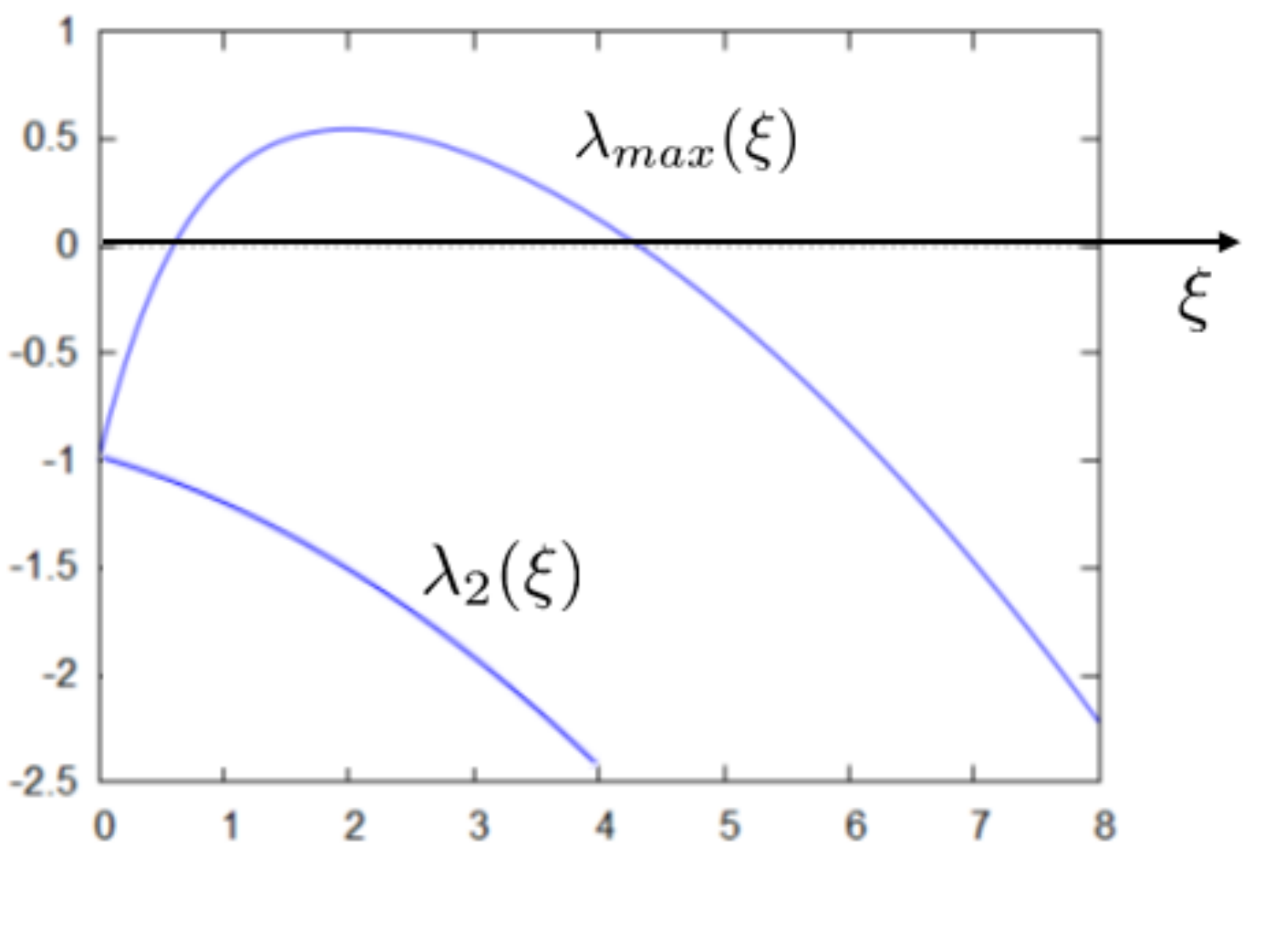}
 \caption{
 Eigenvalues $\lambda_j(\xi )$ of $B(\xi )$. 
 $c_1 = c_2 = 1$, $c_3 = 4$, $c_4 = 3$, $d_1 = 0.05$ and
 $d_2 = 3$ in \rf{RD3-2}. 
 }
 \label{fig7}
\end{figure}

Here we compute the asymptotic profile of $\lambda_{max}(\xi )$
as $|\xi | \rightarrow \infty$. When $|\lambda_{max}(\xi ) 
- \lambda_{\pm\infty}(\xi )|
\rightarrow 0$ as $\xi \rightarrow \pm \infty$, $\lambda_{\pm\infty}(\xi )$
are called the ``asymptotic profiles of $\lambda_{max}(\xi )$''  or the ``asymptotic convergence to
$\lambda_{\pm\infty}(\xi )$'' as $\xi \rightarrow \pm\infty$, which are
expressed by
$\lambda_{max}(\xi ) \rightarrow \lambda_{\pm\infty}(\xi )$
as $\xi \rightarrow \pm\infty$.
In this example,  $\lambda_{max}(\xi )$
has the asymptotic profile
$\lambda_{max}(\xi ) \rightarrow -d_1\xi^2 + c_1$ 
as $|\xi | \rightarrow
\infty$ by virtue of
$d_1 < d_2$ as in Figures \ref{fig2} and \ref{fig7}.
Then we put $\lambda_h (\xi ) 
:= -d_1\xi^2$ (taking the highest order term of $\xi$)
and put $\wh{U} = e^{t\lambda_h(\xi )}\wh{V}$.
We see that $\wh{V}$ satisfies
\begin{equation} \label{RD5}
\wh{V}_t = B_h(\xi )\wh{V},
\end{equation}
where 
$B_h(\xi ) := B(\xi ) - \lambda_h(\xi )I
= 
\left(
\begin{array}{ccc}
c_1 & \; & -c_2 \\
c_3 & \; & -(d_2-d_1)\xi^2 - c_4
\end{array} \right) $
and $I$ is the identity matrix.
We also divide $B_h$ into $B_h(\xi ) = \xi B_1(\xi ) + B_0(\xi )$ with
$B_0(\xi ) = O(1)$ as $|\xi | \rightarrow \infty$, where
$B_1(\xi ) = \left(\begin{array}{ccc}
0 & \; & 0 \\
0 & \; & -(d_2-d_1)\xi
\end{array} \right)$ and 
$B_0(\xi ) = \left(\begin{array}{ccc}
c_1 & \; & -c_2 \\
c_3 & \; & - c_4
\end{array} \right)$ in this example. 
Since the Fourier transformation of the Dirac
$\delta$-function is
$\wh{\delta} (\xi ) = 1$, we approximate it by the heat kernel  $H_{\vep}(x)$
as in Proposition \ref{prop1} and modify $B_h(\xi )$
as $B_{\vep}(\xi ) := \xi B_1(\xi ) + e^{-\vep\xi^2}B_0(\xi )$
for $0 < \vep << 1$ by noting $\wh{H}_{\vep}(\xi ) = e^{-\vep\xi^2}$.
We note that this modification is not necessary
when $|B_0(\xi)| \rightarrow 0$ as $|\xi| \rightarrow \infty$. 
 By using this modified matrix $B_{\vep}(\xi )$, we consider
\begin{equation} \label{RD5-1}
\wh{V}_t = B_{\vep}(\xi )\wh{V}
\end{equation}
instead of \rf{RD5}.

Now introducing small quantity
$0 < \delta << 1$, we consider the equation \rf{RD5-1} at $t+\delta$
\begin{equation} \label{RD6}
\wh{V}_t(t+\delta) = B_{\vep}(\xi )\wh{V}(t+\delta).
\end{equation}
Since $\wh{V}(t) = e^{tB_{\vep}(\xi )}\wh{V}_0$
with the initial data $\wh{V}_0$, $\wh{V}(t+\delta)$ is given by
\[
\wh{V}(t+\delta) 
= e^{\delta B_{\vep}(\xi )}\wh{V}(t) .
\]
Substituting it into \rf{RD6}, we have
\begin{equation} \label{RD6-1}
\wh{V}_t(t+\delta) 
= B_{\vep}(\xi )
e^{\delta B_{\vep}(\xi )}\wh{V}(t) .
\end{equation}

\begin{figure}
 \centering
 \includegraphics[width=7cm, bb=0 0 743 417]{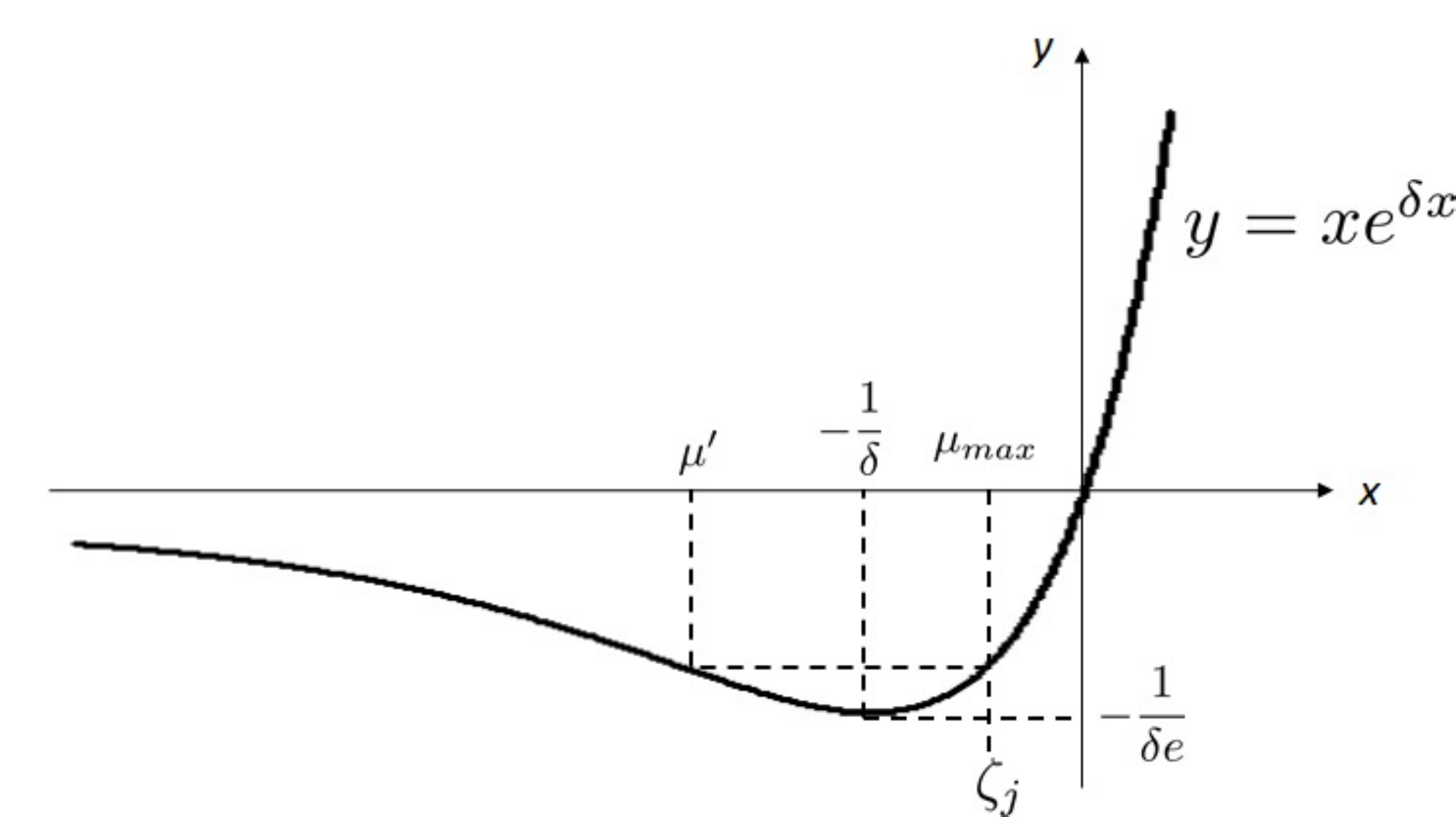}
 \caption{
 Conceptual figure of solutions $\mu$ of \rf{RD6-3}.
 }
 \label{fig6-1-1}\end{figure}

\begin{figure}
 \centering
 \includegraphics[width=12cm,  bb=0 0 953 444]{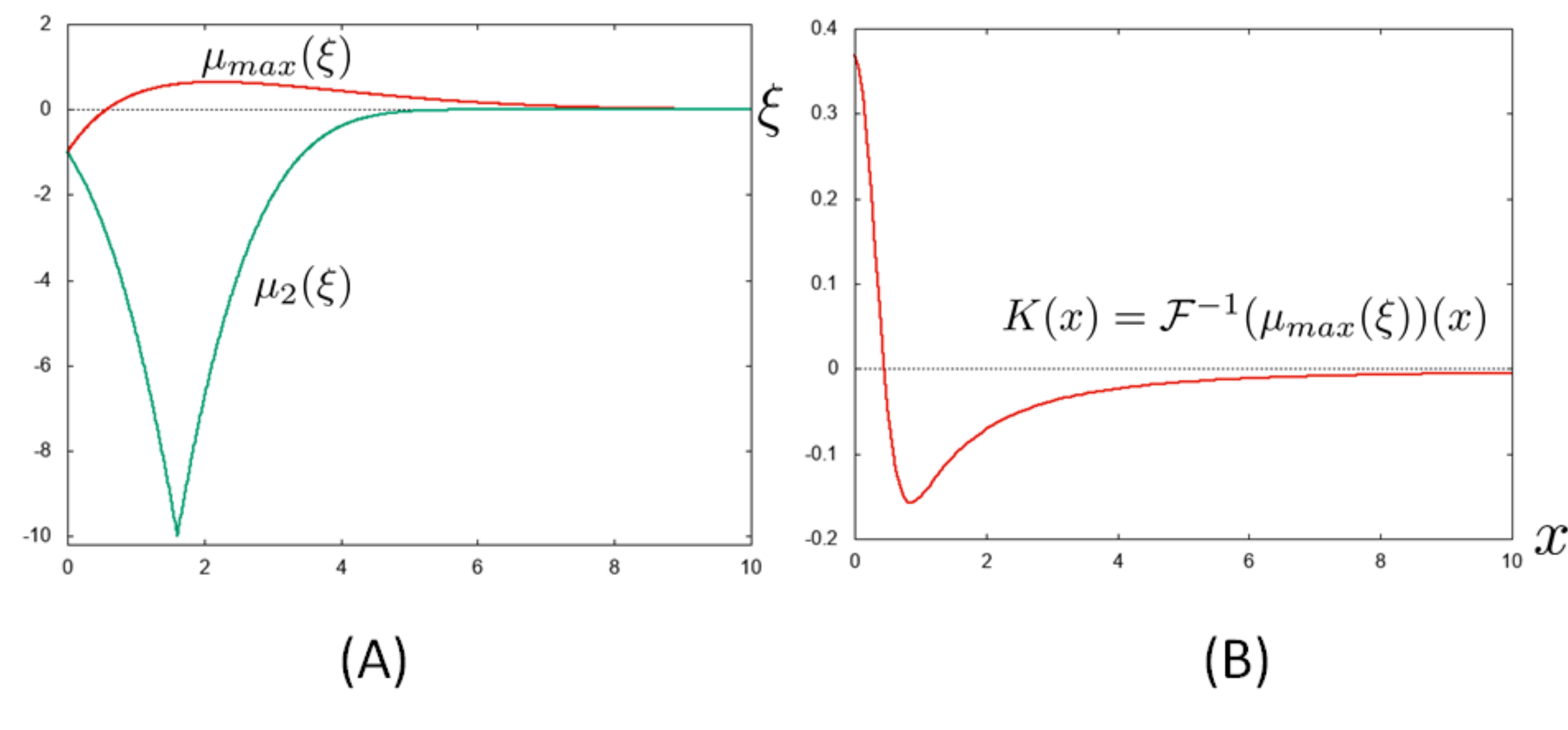}
 \caption{
 (A) Graph of $\mu_{max}(\xi )$. $\mu_2(\xi )$ denotes the
 second largest $Re( \mu )$ in \rf{RD6-3}.
 (B) The derived kernel $K(x) = ({\cal F}^{-1}\mu_{max})(x)$.
 It shows the Mexican hat profile.
 Parameters are same as those of Fig.\ref{fig7}: 
  $c_1 = c_2 = 1$, $c_3 = 4$, $c_4 = 3$, $d_1 = 0.05$ and
 $d_2 = 3$ in \rf{RD3-2}. $\vep$ and $\delta$ in \rf{RD6} are
 taken as
 $\vep = 0.05$, $\delta = 0.1$. 
 Only the right half pictures are drawn both in (A) and (B).
 }
 \label{fig6}
\end{figure}

Let $\zeta_j(\xi )$ and $\Phi_j(\xi )$ ($j =1, 2$)
be eigenvalues and the associated
eigenvectors of $B_{\vep}(\xi )$.
Since $\wh{V}(t,\xi )$ is expressed by the linear combination of
eigenfunctions of $B_{\eps}(\xi )$ as
$\wh{V}(t,\xi) = \dsum_j\alpha_j(\xi)e^{t\mu_j(\xi)}\Phi_j(\xi)$
for some $\alpha_j(\xi)$ and $\mu_j(\xi) \in {\bf C}$, 
the substitution of it into
\rf{RD6-1} gives
\begin{equation} \label{RD6-2}
\dsum_j\alpha_j(\xi)\mu_j(\xi) e^{(t+\delta)\mu_j(\xi)}\Phi_j(\xi)
= \dsum_j\alpha_j(\xi)B_{\vep}(\xi )
e^{(t+\delta )\mu_j(\xi) B_{\vep}(\xi )}\Phi_j (\xi),
\end{equation}
which leads
\[
\mu_j(\xi) e^{\delta\mu_j(\xi)}\Phi_j(\xi )
= \zeta_j(\xi )
e^{\delta\zeta_j(\xi )}\Phi_j(\xi )
\]
and eventually
\begin{equation} \label{RD6-3}
\mu_j(\xi) e^{\delta\mu_j(\xi)}
= \zeta_j(\xi )
e^{\delta\zeta_j(\xi )}.
\end{equation}
Note that $\mu_j(\xi) = \zeta_j(\xi )$
is a trivial solution of \rf{RD6-3} (see e.g. Fig.\ref{fig6-1-1})
and there are infinite many $\mu_j(\xi )$ satisfying \rf{RD6-3}
for any given $\zeta_j(\xi )$.
We need $\mu_j(\xi )$ with maximal real parts among them, that is,
we take $\mu_{max, j}(\xi )$ satisfying $Re(\mu_{max,j}(\xi ))
= \dmax\{ Re(\mu );\; \mu e^{\delta\mu}
= \zeta_j(\xi )e^{\delta\zeta_j(\xi )} \}$, which is known to be real for $\zeta_j(\xi )
\in {\bf R}$ by Proposition\ref{propLW}. Defining 
$\mu_{max}(\xi ) := \dmax_j\{\mu_{max,j}(\xi) ) \}$, we see
$\mu_{max}(\xi )$ is attained by the principal branch of the Lambert W
function as
$\mu_{max}(\xi ) = \frac{1}{\delta}\cdot\dmax_j\{ 
W_0(\delta \zeta_j(\xi )e^{\delta\zeta_j(\xi )})\}$ (see 
Proposition\ref{propLW}).
Since $\mu_{max}(\xi )$ is the most unstable mode,
$\wh{V}(t,\xi)$ of \rf{RD6-1} satisfies
\begin{equation} \label{RD6-4}
|\wh{V}(t,\xi) - \alpha (\xi )e^{t\mu_{max}(\xi )}
\Psi(\xi )| \rightarrow 0 \; (t \rightarrow \infty )
\end{equation}
for almost all initial data $\wh{V}(0,\xi )$
with some $\alpha (\xi )$ and $\Psi (\xi )$.
When \rf{RD6-4} holds, it is called that
$\wh{V}(t,\xi)$ ``asymptotically converges''
$\alpha (\xi )e^{t\mu_{max}(\xi )}\Psi(\xi )$
as $t \rightarrow \infty$.
Then we see
\begin{equation} \label{RD7}
\wh{U}(t,\xi ) =  e^{t\lambda_h(\xi )}\wh{V}(t,\xi )
\rightarrow \alpha (\xi )e^{t(\lambda_h(\xi )+\mu_{max}(\xi ))}
\Psi (\xi )
\end{equation}
as $t \rightarrow \infty$.
Thus, $\wh{U}(t)$ satisfies asymptotically
$\wh{U}_t = (\lambda_h(\xi )+\mu_{max}(\xi ))\wh{U}$
and therefore for any element of $\wh{U}$, say $\wh{w}$, 
$\wh{w}_t = (\lambda_h(\xi )+\mu_{max}(\xi ))\wh{w}$
holds.
Since $\zeta_j(\xi ) \rightarrow 0$ or $-\infty$
as $|\xi | \rightarrow \infty$ by the form of $B_{\vep}(\xi )$,
$\mu_{max}(\xi ) \rightarrow 0$ holds as $|\xi | \rightarrow \infty$
and 
$\mu_{max} \in L^2(\bm{R})$ is expected for almost all cases as shown
in several applications mentioned in this paper. Thus
we obtain the equation of $w$
\begin{equation} \label{RD8}
w_t = d_1w_{xx} + K*w
\end{equation}
as the effective equation,
where we use ${\cal F}^{-1}(-d_1\xi^2) = d_1\pa_x^2$ and
define
$K(x) := {\cal F}^{-1}(\mu_{max}(\xi ))(x)$,
the inverse Fourier transformation
of the function $\mu_{max}(\xi )$.
Fig.\ref{fig6} is a numerical simulation under the indicated parameters.
We note that the kernel $K(x)$
gives the Mexican hat type profile.
In practical computations, we use
\begin{equation} \label{RD8-1}
w_t = d_1w_{xx} + \chi (K*w)
\end{equation}
with an appropriate cutt-off function $\chi (r)$ as in \rf{KT1}.

\begin{REM}
The above matrix $B(\xi )$ and the eigenvalue $\lambda_j(\xi )$ can be written as 
$B(\xi^2 )$ and $\lambda_j(\xi^2 )$ from the forms.
In general,  
the matrix and the eigenvalue corresponding to $B(\xi )$ and $\lambda_j(\xi )$
depend only on $\xi^2$
when the spatial movement is only by a local diffusion like this example.
In this case,
one dimensional analysis is enough because the Fourier transformation
on 2D just replaces $\xi^2$ by $\xi^2+\eta^2$, that is, the matrix $B(\xi^2+\eta^2 )$
and the eigenvalue $\lambda_j(\xi^2+\eta^2 )$ in 2D
without any other change, which means the reduced kernel
in 1D and 2D are almost same (refer Section 5).
\end{REM}
\begin{REM}
The equations \rf{KT1} and \rf{RD8-1} seem different
each other. But they are essentially same in the following sense:
The diffusion is approximated by the relation
$J_{\eps}*u = \gamma_0u + \eps^2\gamma_1u_{xx}  + O(\eps^4)$ as
$\eps \rightarrow +0$, where
$J_{\eps}(x ) := \frac{1}{\eps}J(\frac{x }{\eps})$
for an even function $J(x)$ satisfying $J(x) \ge 0$ and
$\gamma_0 := \dint_{\bf R}J(x)dx$,
$\gamma_1 := \frac{1}{2}\dint_{\bf R}x^2J(x)dx$ (e.g. \cite{BCC, Mu}).
Then the reduced equation \rf{RD8} without cut-off function in \rf{RD8-1}
is expressed as
\begin{equation} \label{RD8-2}
w_t = K'*w -\frac{\gamma_0}{\eps^2\gamma_1}w + O(\eps^2)
\end{equation}
with $K'(x) := \frac{d_1}{\eps^2\gamma_1}J_{\eps}(x) + K(x)$, 
which is the same form as \rf{KT1} without cut-off function.
Thus, \rf{KT1} and \rf{RD8-1} are essentially same and we
are adopting the equations in the form of \rf{RD8-1} in this paper.  
\end{REM}

\SECT{Treatment for arbitrary networks}

\subsection{Basic treatment for general networks \\
- case of real eigenvalues -}

In this section, we first give the idea to treat general cases of multiple components, 
but
consider on one dimensional space just for the simplicity.
Let $U = {^t}(u_1, \cdots , u_N) \in \bm{R}^N$ and $U$
is supposed to be described by the equation
\begin{equation} \label{e1}
U_t = {\bm J}U + AU,
\end{equation}
where ${\bm J}U := {^t}(J_1*u_1, \cdots , J_N*u_N)$ for
${\bm J} := {^t}(J_1, \cdots , J_N)$ with kernels $J_j$ corresponds
to the spatial movement and $N$-th order square matrix $A$
denotes the local network  like the matrix $A$ in \rf{RD4}.
Taking the Fourier transformation of \rf{e1}, we have
\begin{equation} \label{e1-1}
\wh{U}_t = B(\xi )\wh{U},
\end{equation}
where $B(\xi ) := \wh{\bm{J}} + A$. Let $\lambda_j(\xi )$
($j = 1, \cdots, N$)
be eigenvalues of $B(\xi )$.
In this subsection, all $\lambda_j(\xi )$ are 
 supposed to be real.
Quite similarly to the previous section, we define
$\lambda_{max}(\xi ) := \dmax_j\{\lambda_j(\xi )\}$
and $\lambda_h(\xi )$ be the highest order term in $\xi$
of the asymptotic profile of
$\lambda_{max}(\xi )$ as $|\xi | \rightarrow \infty$.
Here we give the precise way to get $\lambda_h(\xi )$.
Let $\lambda_{\pm\infty}(\xi )$ be the asymptotic profile of 
$\lambda_{max}(\xi )$ satisfying $\lambda_{max}(\xi ) \rightarrow
\lambda_{\pm\infty}(\xi )$ as $\xi \rightarrow \pm\infty$.
When $|\lambda_{\pm\infty}(\xi )| 
\rightarrow \infty$ as $\xi \rightarrow \pm\infty$ and
the highest order terms of $\lambda_{\pm\infty}(\xi )$
are coincident each other, 
take the same highest
order term of $\lambda_{\pm\infty}(\xi )$ as $\lambda_h(\xi )$ and
when $\lambda_{\pm\infty}(\xi ) \rightarrow O(1)$, take
$\lambda_h(\xi ) = 0$.
For example, when $\lambda_{\pm\infty}(\xi )$ is given by
$\lambda_{\pm\infty}(\xi ) = a_1\xi^2 + a_{\pm 2}\xi + a_{\pm 3}$ 
for $a_j \in {\bf C}$
and $a_1 \ne 0$,
we put $\lambda_h(\xi ) = a_1\xi^2$ and
when $\lambda_{\pm\infty}(\xi ) = a_{\pm 3}$, 
we put
$\lambda_h(\xi ) = 0$.

Transforming $\wh{U}(t) 
= e^{t\lambda_h(\xi )}\wh{V}(t)$,
we have the equation of $\wh{V}$
\begin{equation} \label{e2}
\wh{V}_t = B_h(\xi )\wh{V}
\end{equation}
with $B_h(\xi ) := \{ B(\xi ) - \lambda_h(\xi )I\}$.
Expressing $B_h(\xi )$ in the form $B_h(\xi ) = \xi B_1(\xi ) + B_0(\xi )$
with  $B_0(\xi ) = O(1)$
as $|\xi | \rightarrow \infty$, we define 
$B_{\vep}(\xi ) := \xi B_1(\xi ) + e^{-\vep\xi^2}B_0(\xi )$
and consider the equation
\begin{equation} \label{e2-1}
\wh{V}_t = B_{\vep}(\xi )\wh{V}
\end{equation}
instead of \rf{e2}.
Let $\zeta_j(\xi )$ and $\Phi_j(\xi )$ ($j =1,\cdots,N$)
be eigenvalues and the
associated eigenvectors of $B_{\vep}(\xi )$, which are supposed
to be real in this subsection.

Taking the time $t+\delta$ for $0 < \delta << 1$ in \rf{e2-1}, we see
\begin{equation} \label{e3}
\wh{V}_t(t+\delta) = B_{\vep}(\xi )\wh{V}(t+\delta)
\end{equation}
and therefore 
\begin{equation} \label{e3-1}
\wh{V}_t(t+\delta) 
= B_{\vep}(\xi )
e^{\delta B_{\vep}(\xi )}\wh{V}(t)
\end{equation}
from \rf{e2-1}. Then 
we express $\wh{V}(t,\xi ) = \dsum_j\alpha_j(\xi)e^{t\mu_j(\xi)}\Phi_j(\xi)$
and the substitution of it into \rf{e3-1} leads
\begin{equation} \label{e3-3}
\mu_j(\xi) e^{\delta\mu_j(\xi )}
= \zeta_j(\xi )
e^{\delta\zeta_j(\xi )}
\end{equation}
in quite a similar manner to \rf{RD6-2} and \rf{RD6-3}.
Define 
$\mu_{max}(\xi )$ as the value satifying
\[
Re(\mu_{max}(\xi )) = \dmax_j\{Re (\mu);\;
 \mu e^{\delta\mu}= \zeta_j(\xi )e^{\delta\zeta_j(\xi )},\;
 j = 1,2,\cdots, N \}
 \]
 similarly to Section 3,
which is given by $\mu_{max}(\xi ) = \dmax_j\{
1/\delta\cdot W_0(\delta\zeta_j(\xi )e^{\delta\zeta_j(\xi )})\}$ 
when
$\zeta_j(\xi ) \in {\bf R}$ by Proposition \ref{propLW}.
\begin{REM}
$\mu_j(\xi ) 
= \zeta_j(\xi )$ holds
when $\zeta_j(\xi ) > -\frac{1}{\delta}$
(see Figure \ref{fig6-1-1}). 
\end{REM}
Since $\mu_{max}(\xi )$ is the most unstable eigenvalue for $\xi \in {\bf R}$, 
we can expect
$\wh{V}(t)$ of \rf{e3-1} asymptotically converges
$\wh{V}(t,\xi ) \rightarrow \alpha (\xi )e^{t\mu_{max}(\xi )}
\Phi(\xi )$
($t \rightarrow \infty$) with some $\alpha (\xi )$ and $\Phi (\xi )$
for almost all initial data
and hence
\begin{equation} \label{e3-4}
\wh{U}(t,\xi ) =  e^{t\lambda_h(\xi )}\wh{V}(t,\xi )
\rightarrow \alpha (\xi )e^{t(\lambda_h(\xi )+\mu_{max}(\xi ))}
\Phi (\xi )
\end{equation}
as $t \rightarrow \infty$.
Thus, $\wh{U}(t)$ satisfies asymptotically
$\wh{U}_t = (\lambda_h(\xi )+\mu_{max}(\xi ))\wh{U}$
and therefore for any element of $\wh{U}$, say $\wh{w}$, 
$\wh{w}_t = (\lambda_h(\xi )+\mu_{max}(\xi ))\wh{w}$
holds. Consequently 
we obtain the equation of $w$
\begin{equation} \label{e3-5}
w_t = {\cal L}w + K*w
\end{equation}
as the effective equation,
where  $K(x) := {\cal F}^{-1}(\mu_{max}(\xi ))$ and
${\cal L} := {\cal F}^{-1}(\lambda_h(\xi ))$, which
is treated as a differential operator such as
${\cal F}^{-1}(-\xi^2 ) = \pa_x^2$.
Here we note that $\mu_{max} (\xi ) \rightarrow 0$ as $|\xi | \rightarrow \infty$
because $\zeta_j(\xi ) \rightarrow 0$ or $-\infty$.

This is the basic idea in the case that all $\lambda_j(\xi )$ are 
real valued and  there exists
only one maximal eigenvalue for each $\xi$.
But this is a specially simple case because there can be several other cases.
In the following subsection, the case that
$B(\xi )$ have complex eigenvalues 
on some range of $\xi \in {\bf R}$ is considered.

\subsection{Basic treatment for general networks \\
- case of complex eigenvalues -}
In this section, we consider the case that $B(\xi )$ in \rf{e1-1} have
complex eigenvalues for an interval of $\xi $.
We define $\lambda_{max}(\xi ) := \dmax_j\{ Re(\lambda_j(\xi )) \}$ and
the asymptotic profile $\lambda_{\infty}(\xi )$ of $\lambda_{max}(\xi )$
by $\lambda_{max}(\xi ) \rightarrow \lambda_{\infty}(\xi )$ as $|\xi | \rightarrow
\infty$. First, we assume 
\begin{equation} \label{e10-1-2}
\lambda_j(\xi ) - \lambda_{\infty}(\xi ) \rightarrow 0\; (|\xi | \rightarrow \infty ).
\end{equation}
Let $\lambda_h(\xi )$ be the highest order term of $\lambda_{\infty}(\xi )$ and
$B_h(\xi ) := B(\xi ) - \lambda_h(\xi )I$. 
As in the previous section, we divide $B_h(\xi )$ as
$B_h(\xi ) = \xi B_1(\xi ) + B_0(\xi )$ with $B_0(\xi ) = O(1)$ and define
$B_{\vep}(\xi ) = \xi B_1(\xi ) + e^{-\vep\xi^2}B_0(\xi )$.

Let $\zeta_j(\xi )$ and $\Phi_j(\xi )$ be eigenvalues and the associated 
eigenvectors of $B_{\vep}(\xi )$
and consider the equation \rf{e3-1} together with \rf{e3-3}
in the case of complex eigenvalues. Define 
\[
\mu_{max}(\xi ) := \dmax_j\{Re (\mu);\;
 \mu e^{\delta\mu}= \zeta_j(\xi )e^{\delta\zeta_j(\xi )},\;
 j = 1,2,\cdots, N \}
\]
and assume that
$\zeta_1(\xi)$, $\zeta_2(\xi)$ 
satisfy
the following three assumptions (H1), (H2) and (H3). 

\ \\
(H1): $\mu_{max}(\xi )$ is attained by $\zeta_1(\xi )$, that is,
$\mu_{max}(\xi ) = \max\{ Re(\mu );\;
\mu e^{\delta\mu} = \zeta_1(\xi )e^{\delta\zeta_1(\xi )}\}$.
\\
(H2):  There exists $\xi_c$ such that $\zeta_1(\xi ) > \zeta_2(\xi )$
for $\xi < \xi_c$, $\zeta_1(\xi_c) = \zeta_2(\xi_c) \in {\bf R}$ and
$\zeta_1(\xi ) = \ovl{\zeta_2(\xi )} = a(\xi ) + ib(\xi )$for $\xi > \xi_c$
with $a(\xi ) \in {\bf R}$ and $b(\xi ) >  0$ satisfying
$a(\xi_c) = \zeta_1(\xi_c) = \zeta_2(\xi_c)$ and $b(\xi_c) =0$.
Moreover, $\Phi_1(\xi_c) = \Phi_2(\xi_c)$ holds.
\\
(H3): 
Define
$\Psi (\xi ) := \ddfrac{\Phi_2(\xi ) - \Phi_1(\xi )}{\zeta_2(\xi ) - \zeta_1(\xi )}$. Then, there exists a limit $\dlim_{\xi\rightarrow\xi_c}\Psi (\xi )$.\\
\begin{REM}
It follows  that
$b(\xi ) \rightarrow 0$ as $|\xi | \rightarrow \infty$ by \rf{e10-1-2}
and also that
$\Psi (\xi )$ is a continuous real vector-valued function 
on the whole line ${\bf R}$ by the definition in (H3).
\end{REM}
Substituting
$\Phi_2(\xi ) = \Phi_1(\xi ) + \{\zeta_2(\xi ) - \zeta_1(\xi )\}\Psi (\xi )$
into $B_{\vep}(\xi )\Phi_2(\xi ) 
= \zeta_2(\xi )\Phi_2(\xi )$ for $\xi < \xi_c$, we
have 
$\zeta_1(\xi )\Phi_1(\xi ) + \{ \zeta_2(\xi ) - \zeta_1(\xi  )\}
B_{\vep}(\xi )\Psi (\xi )
= \zeta_2(\xi )\Phi_1(\xi ) 
+ \zeta_2(\xi )\{\zeta_2(\xi ) - \zeta_1(\xi  )\} B_{\vep}(\xi )\Psi (\xi )$
and consequently $\{ B_{\vep}(\xi ) - \zeta_2(\xi )I\}\wt{\Phi}(\xi ) = \Phi_1(\xi )$.
That is, 
\begin{equation} \label{e10-1}
\left\{ \begin{array}{lcl}
\{ B_{\vep}(\xi ) - \zeta_1(\xi )I\}\Phi_1 (\xi ) &=& 0, \\
\{ B_{\vep}(\xi ) - \zeta_2(\xi )I\}\Psi (\xi ) &=& \Phi_1(\xi )
\end{array} \right.
\end{equation}
holds for $\xi < \xi_c$.

On the other hand, we have for $\xi > \xi_c$,
\[
\left\{ \begin{array}{lcl}
\{ B_{\vep}(\xi ) - a(\xi )I\}{\bf p}(\xi ) &=& -b(\xi ){\bf q}(\xi ), \\
\{ B_{\vep}(\xi ) - a(\xi )I\}{\bf q}(\xi ) &=& b(\xi ){\bf p}(\xi )
\end{array} \right.
\]
when we express $\Phi_1(\xi ) = {\bf p}(\xi ) + i{\bf q}(\xi ) 
= \ovl{\Phi_2(\xi )}$ with ${\bf p}(\xi ),\; {\bf q}(\xi ) \in {\bf R}^N$.
Since the definition of $\Psi (\xi )$ leads
${\bf q}(\xi ) = b(\xi )\Psi (\xi )$ and consequently
\begin{equation} \label{e10-2}
\left\{ \begin{array}{lcl}
\{ B_{\vep}(\xi ) - a(\xi )I\}{\bf p}(\xi ) &=& -b^2(\xi )\Psi (\xi ), \\
\{ B_{\vep}(\xi ) - a(\xi )I\}\Psi (\xi ) &=& {\bf p}(\xi )
\end{array} \right.
\end{equation}
for $\xi > \xi_c$. 
\rf{e10-1} and \rf{e10-2} imply that
the functions $b(\xi )$, $\zeta_j(\xi )$ ($j=1,2$) and the vector
$\Phi_1(\xi )$ can be extended
as real valued continuous functions and a vector on the whole line ${\bf R}$ by 
\begin{equation} \label{e10-21}
b(\xi ) := \left\{ \begin{array}{cc} 
0 & (\xi \le \xi_c) \\
b(\xi ) & (\xi > \xi_c ) \end{array} \right. ,\;
\omega_j (\xi ) := \left\{ \begin{array}{cc} 
\zeta_j(\xi ) & (\xi \le \xi_c) \\
a(\xi ) & (\xi > \xi_c ) \end{array} \right. ,\;
\Psi_1(\xi ) := \left\{ \begin{array}{cc} 
\Phi_1(\xi ) & (\xi \le \xi_c) \\
{\bf p}(\xi ) & (\xi > \xi_c ) \end{array} \right. .
\end{equation}
Thus \rf{e10-1} and \rf{e10-2} are unified as one system on ${\bf R}$
\begin{equation} \label{e10-3}
\left\{ \begin{array}{lcl}
\{ B_{\vep}(\xi ) - \omega_1 (\xi )I\}\Psi_1(\xi ) &=& -b^2(\xi )\Psi , \\
\{ B_{\vep}(\xi ) - \omega_2(\xi )I\}\Psi (\xi ) &=& \Psi_1(\xi ) .
\end{array} \right.
\end{equation}

Now we go back the equation \rf{e3-1} of $\wh{V}$.
Let
the solution $\wh{V}$ of \rf{e3-1} be of  the form
$\alpha \Psi_1(\xi ) + \beta\Psi(\xi )$
for $\alpha,\; \beta \in {\bf R}$. Since
\rf{e10-3} shows
\[
B_{\vep}(\xi )(\alpha \Psi_1(\xi ) + \beta\Psi (\xi ))
= (\omega_1(\xi )\alpha + \beta )\Psi_1(\xi )
+ (-b^2(\xi )\alpha + \omega_2(\xi )\beta )\Psi (\xi ) ,
\]
\rf{e3-1} becomes the equation of $\alpha$ and
$\beta$ as
\begin{equation} \label{e10-4}
\ode{}{t}\left( \begin{array}{c}
\alpha(t+\delta ) \\
\beta(t+\delta) 
\end{array} \right)
= \wt{B}_{\vep}(\xi )e^{\delta\wt{B}_{\vep}(\xi )}
\left( \begin{array}{c}
\alpha(t) \\
\beta(t) 
\end{array} \right),
\end{equation}
where $\wt{B}_{\vep}(\xi ) := \left( \begin{array}{cc}
\omega_1(\xi ) & 1 \\
-b^2(\xi ) & \omega_2(\xi ) \end{array} \right)$.
Since the eigenvalues of $\wt{B}_{\vep}(\xi )$ are
$\wt{\omega}_1(\xi ) := \omega_1(\xi ) + ib(\xi )$
and $\wt{\omega}_2(\xi ) := \omega_2(\xi ) - ib(\xi )$
with associated eigenvectors $\wt{\Psi}_1(\xi ) =
{^t}(1, ib(\xi ))$ and 
$\wt{\Psi}_2(\xi ) = {^t}(1, \omega_2(\xi ) - \omega_1(\xi ) -ib(\xi ))$
respectively, 
$\left(
\begin{array}{c}
\alpha \\ \beta \end{array}\right) = e^{t\mu}\wt{\Psi}$ with \rf{e10-4}
leads
\begin{equation} \label{e10-41}
\mu e^{\delta\mu} = \wt{\omega}_je^{\delta\wt{\omega}_j}
= e^{\delta\omega_j}\{
(\omega_j\cos\delta b - b\sin\delta b ) \pm
i(b\cos\delta b + \omega_j\sin\delta b)\}
\end{equation}
by taking $\wt{\Psi} = \wt{\Psi}_j(\xi )$.
Let the set $\wt{{\cal M}}_j (\xi ) := \{ \mu ;\;
\mu e^{\delta\mu} = \wt{\omega}_je^{\delta\wt{\omega}_j}\}$
and define $\wt{\mu}_j(\xi ) \in \wt{{\cal M}}_j(\xi )$ for $j = 1, 2$
such that $Re(\wt{\mu}_j(\xi )) 
= \max\{ Re(\mu );\; \mu \in \wt{{\cal M}}_j(\xi )\}$.
Here we note that these $\wt{\mu}_j(\xi )$ give $\mu_{max}(\xi )$,
that is, $\mu_{max}(\xi ) = \dmax_{j=1,2}\{ Re(\wt{\mu}_j(\xi ))\}$.
Then we can take the solution $\left(
\begin{array}{c}
\alpha \\ \beta \end{array}\right)$ of \rf{e10-4} by
$\left(
\begin{array}{c}
\alpha (t) \\ \beta (t) \end{array}\right) =
e^{t\wt{\mu}_1}\wt{\Psi}_1 + e^{t\wt{\mu}_2}
\wt{\Psi}_2$, which
leads the equation of $\alpha (t)$ and $\beta (t)$ as
\begin{eqnarray} \label{e10-42}
\ode{}{t}\left( \begin{array}{c}
\alpha \\
\beta
\end{array}
 \right) &=& [\wt{\Psi}_1,\wt{\Psi}_2]
 \left(\begin{array}{cc} \wt{\mu}_1 & 0 \\
 0 & \wt{\mu}_2 \end{array}\right)
  [\wt{\Psi}_1,\wt{\Psi}_2]^{-1}
 \left( \begin{array}{c}
\alpha \\
\beta
\end{array}
 \right)  \\
 \nonumber
 &=&
 \ddfrac{1}{\omega_2 - \omega_1 - 2ib}
 \left( \begin{array}{cc}
 (\omega_2-\omega_1)\wt{\mu}_1 - ib(\wt{\mu}_1+\wt{\mu}_2) &
 -\wt{\mu}_1 + \wt{\mu}_2 \\
 ib(\omega_2 - \omega_1 - ib)(\wt{\mu}_1-\wt{\mu}_2) &
 (\omega_2-\omega_1)\wt{\mu}_2 - ib(\wt{\mu}_1+\wt{\mu}_2)
 \end{array}\right)
 \left( \begin{array}{c}
\alpha \\
\beta
\end{array}
 \right) .
 \end{eqnarray}

 By the way of construction in \rf{e10-21}, we see
 $\wt{\omega}_1(\xi) \ne \wt{\omega}_2(\xi)$
 and $\wt{\omega}_j(\xi) = \zeta_j(\xi ) \in {\bf R}$ for $\xi < \xi_c$,
 $\wt{\omega}_1(\xi_c) = \wt{\omega}_2(\xi_c)$,
$\wt{\omega}_1(\xi) = \ovl{\wt{\omega}_2(\xi)}$
for $\xi > \xi_c$,
 we can also assume 
 $\wt{\mu}_1(\xi) \ne \wt{\mu}_2(\xi)$
 and $\wt{\mu}_j(\xi) \in {\bf R}$ for $\xi < \xi_c$,
 $\wt{\mu}_1(\xi_c) = \wt{\mu}_2(\xi_c)$,
$\wt{\mu}_1(\xi) = \ovl{\wt{\mu}_2(\xi)}$
for $\xi > \xi_c$.
 Then, 
 \rf{e10-42} becomes
 \begin{equation} \label{e10-43}
 \ode{}{t}\left( \begin{array}{c}
\alpha \\
\beta
\end{array}
 \right)
 = \ddfrac{1}{\zeta_2 - \zeta_1}
 \left( \begin{array}{cc}
 (\zeta_2-\zeta_1)\wt{\mu}_1 & -\wt{\mu}_1+\wt{\mu}_2 \\
 0 &  (\zeta_2-\zeta_1)\wt{\mu}_2
 \end{array} \right)\left( \begin{array}{c}
\alpha \\
\beta
\end{array}
 \right)
 = 
 \left( \begin{array}{cc}
 \wt{\mu}_1 & \frac{\wt{\mu}_2-\wt{\mu}_1}{\omega_2-\omega_1} \\
 0 &  \wt{\mu}_2
 \end{array} \right)\left( \begin{array}{c}
\alpha \\
\beta
\end{array}
 \right)
 \end{equation}
 for $\xi < \xi_c$ and
 \begin{equation} \label{e10-44}
 \ode{}{t}\left( \begin{array}{c}
\alpha \\
\beta
\end{array}
 \right)
 = \ddfrac{1}{-2ib}
 \left( \begin{array}{cc}
 -ib(\wt{\mu}_1+\wt{\mu}_2) & -\wt{\mu}_1+\wt{\mu}_2 \\
 b^2(\wt{\mu}_1-\wt{\mu}_2) &  -ib(\wt{\mu}_1+\wt{\mu}_2)
 \end{array} \right)
 \left( \begin{array}{c}
\alpha \\
\beta
\end{array}
 \right)
 =
 \left( \begin{array}{cc}
 \wt{a} & \wt{b}/b \\
 -\wt{b}b & \wt{a}
 \end{array} \right)
 \left( \begin{array}{c}
\alpha \\
\beta
\end{array}
 \right)
 \end{equation}
 for $\xi > \xi_c$, where $\wt{\mu}_1(\xi ) = \ovl{\wt{\mu}_2(\xi )} 
 =\wt{a}(\xi )+i\wt{b}(\xi )$ for $\xi > \xi_c$.
 Therefore, if the limit exists and
\[
\dlim_{\xi\rightarrow\xi_c-0}
\frac{\wt{\mu}_2(\xi )-\wt{\mu}_1(\xi )}{\zeta_2(\xi )-\zeta_1(\xi )}
=
\dlim_{\xi\rightarrow\xi_c+0}\frac{\wt{b}(\xi )}{b(\xi )} =: \wt{p}_c
\] 
holds, then we can define real valued continuous functions on ${\bf R}$ by
\begin{equation} \label{e10-45}
\wt{p}(\xi ) := \left\{ \begin{array}{cc} 
\frac{\wt{\mu}_2(\xi )-\wt{\mu}_1(\xi )}{\zeta_2(\xi )-\zeta_1(\xi )}
 & (\xi < \xi_c) \\
 \wt{p}_c & (\xi = \xi_c) \\
\frac{\wt{b}(\xi )}{b(\xi )} & (\xi > \xi_c ) \end{array} \right. ,\;
\wt{\nu}_j (\xi ) := \left\{ \begin{array}{cc} 
\wt{\mu}_j(\xi ) & (\xi \le \xi_c) \\
\wt{a}(\xi ) & (\xi > \xi_c ) \end{array} \right. ,\;
\wt{q}(\xi ) := \left\{ \begin{array}{cc} 
0 & (\xi \le \xi_c) \\
-\wt{b}(\xi )b(\xi ) & (\xi > \xi_c ) \end{array} \right. 
\end{equation}
and equations \rf{e10-43} and \rf{e10-44} are unified as one equation
\begin{equation} \label{e10-46}
 \ode{}{t}\left( \begin{array}{c}
\alpha \\
\beta
\end{array}
 \right)
 = 
 \left( \begin{array}{cc}
 \wt{\nu}_1(\xi ) & \wt{p}(\xi ) \\
 \wt{q}(\xi ) & \wt{\nu}_2(\xi )
 \end{array} \right)
 \left( \begin{array}{c}
\alpha \\
\beta
\end{array}
 \right) .
\end{equation}
Thus we get $\wh{V}$ of \rf{e3-1} as
$\wh{V}= \alpha\Psi_1 + \beta\Psi$ and
$\wh{U} = e^{t\lambda_h}\wh{V} 
= e^{t\lambda_h}\alpha\Psi_1 + e^{t\lambda_h}\beta\Psi$.
This implies that $\hat{U} = \alpha'\Psi_1 + \beta'\Psi$ satisfies
\begin{equation} \label{e10-47}
 \ode{}{t}\left( \begin{array}{c}
\alpha' \\
\beta'
\end{array}
 \right)
 = 
 \left( \begin{array}{cc}
 \wt{\nu}_1(\xi ) + \lambda_h(\xi ) & \wt{p}(\xi ) \\
 \wt{q}(\xi ) & \wt{\nu}_2(\xi ) + \lambda_h(\xi )
 \end{array} \right)
 \left( \begin{array}{c}
\alpha' \\
\beta'
\end{array}
 \right) .
\end{equation}

Consequently, we can get the equation of $X := {\cal F}^{-1}\alpha'$,
$Y := {\cal F}^{-1}\beta'$ as
\begin{equation} \label{e10-5}
\left\{ \begin{array}{ccl}
\dot{X} &=& {\cal L}X + K * X + L*Y, \\
\dot{Y} &=& {\cal L}Y + M * X + N * Y,
\end{array} \right.
\end{equation}
where $K(x) := {\cal F}^{-1}(\wt{\nu}_1 (\xi ) )(x)$, 
$L(x) := {\cal F}^{-1}(\wt{p}(\xi ))(x)$,
$M(x) := {\cal F}^{-1}(\wt{q}(\xi ))(x)$,
$N(x) := {\cal F}^{-1}(\wt{\nu}_2(\xi ))(x)$
and the differential operator
${\cal L} := {\cal F}^{-1}(\lambda_h(\xi ))$.

The simpler case stated in the previous subsection
corresponds to the case $Y = 0$ and $M = 0$ in \rf{e10-5}. 


We summarize several cases which can occur in practical situations
in the following subsections.

\subsection{Practical ways}

In practical situations, 
we will meet several difficulties and complication in calculation.
In this section, we give practical ways to obtain effective equations
by using approximations.

Hereafter, we assume the asymptotic profile $\lambda_{\infty}(\xi )$
is real for any $\xi \in {\bf R}$.

\begin{figure}
 \centering
 \includegraphics[width=12cm,  bb=0 0 928 414]{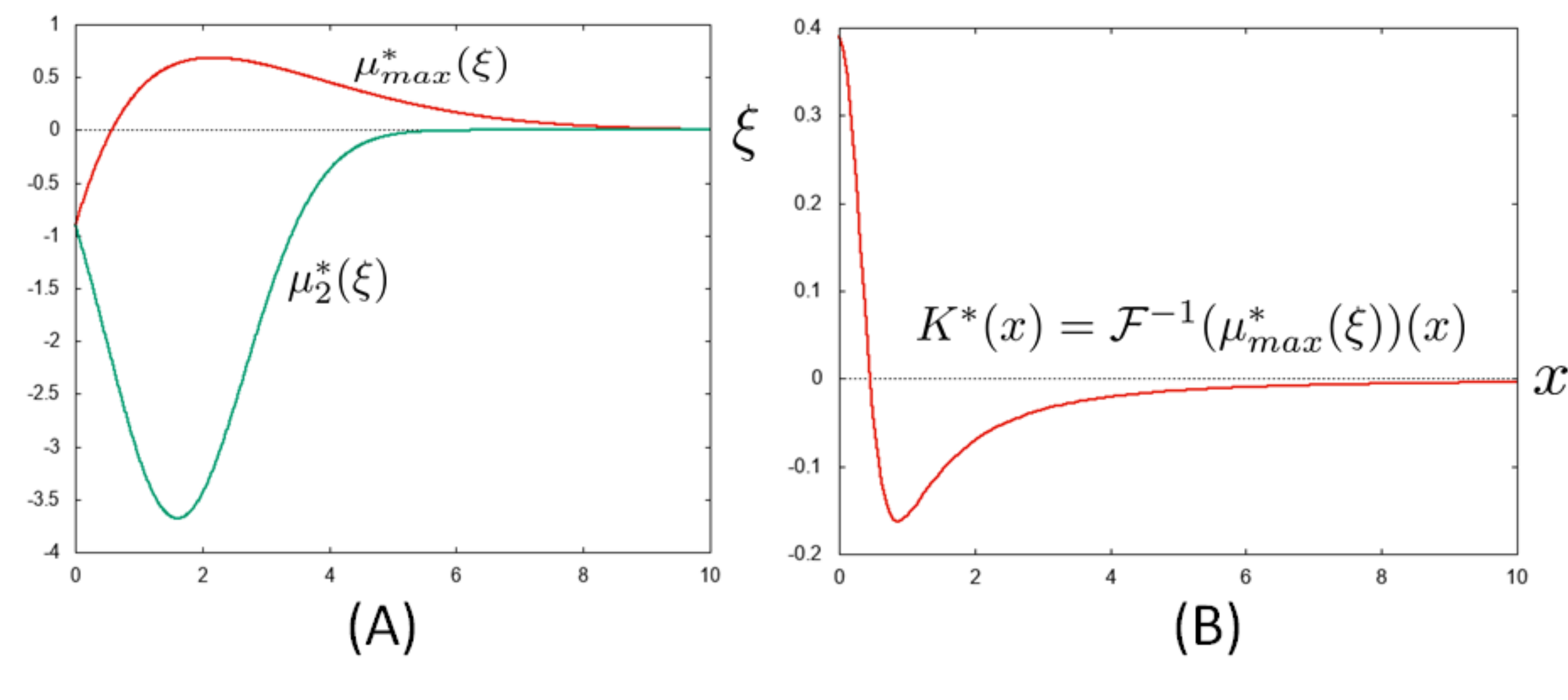}
 \caption{
 (A) Graph of $\mu^*_{max}(\xi )$. $\mu^*_2(\xi )$ denotes the
 second largest $M^*(\zeta_j(\xi ))$ for the eigenvalues $\zeta_j(\xi )
 \in B_{\vep}(\xi )$.
 (B) The derived kernel $K^*(x) = ({\cal F}^{-1}\mu^*_{max})(x)$,
 which shows the Mexican hat profile.
 Parameters are same as those of Fig.\ref{fig7} and Fig.\ref{fig6}: 
  $c_1 = c_2 = 1$, $c_3 = 4$, $c_4 = 3$, $d_1 = 0.05$,
 $d_2 = 3$ and  
 $\vep = 0.05$, $\delta = 0.1$. 
 Only the right half pictures are drawn both in (A) and (B).
 }
 \label{fig6-1}
\end{figure}

In order to  obtain $\mu_{max}(\xi )$ in \rf{e3-3}, \rf{e10-41},
we give the following two practical ways.

\ \\
Practical way I):\\
Define $M_{max}(z) := \max\{ Re(\mu );\; \mu e^{\delta\mu} = ze^{\delta z} \}$ for a given $z \in {\bf C}$.
Then $\mu_{max}(\xi )$ is given by
$\mu_{max}(\xi ) = \dmax_j\{ M_{max}(\zeta_j(\xi ));\;
\zeta_j(\xi ) \in \sigma (B_{\vep}(\xi ))\}$.
Thus, $\mu_{max}(\xi )$ is calculated by the function $M_{max}(z)$.
On the other hand, $M_{max}(z) = \frac{1}{\delta}W_0(\delta ze^{\delta z})$ holds
for $z \in {\bf R}$ by Proposition \ref{propLW} and the graph of
$\frac{1}{\delta}W_0(\delta ze^{\delta z})$
is qualitatively
close to $M^*(z) := ze^{\delta z}$ when $z \in {\bf R}$.
Hence we adopt $M^*(z)$ instead of $M_{max}(z)$ as a rough approximation
and use $\mu^*_{max}(\xi ) := \dmax_j\{ M^*(\zeta_j(\xi ))\}$
as the approximation of $\mu_{max}(\xi )$
when $\zeta_j(\xi ) \in {\bf R}$. 
In Fig.\ref{fig6-1}, the kernel
$K^*(x) = {\cal F}^{-1}(\mu^*_{max}(\xi ))(x)$ is drawn under
the same parameters as Fig.\ref{fig6}. It suggests that 
the essential structure of the original kernel $K(x)$ is retained by this rough
approximation.

\begin{figure}
 \centering
 \includegraphics[width=12cm,  bb=0 0 903 410]{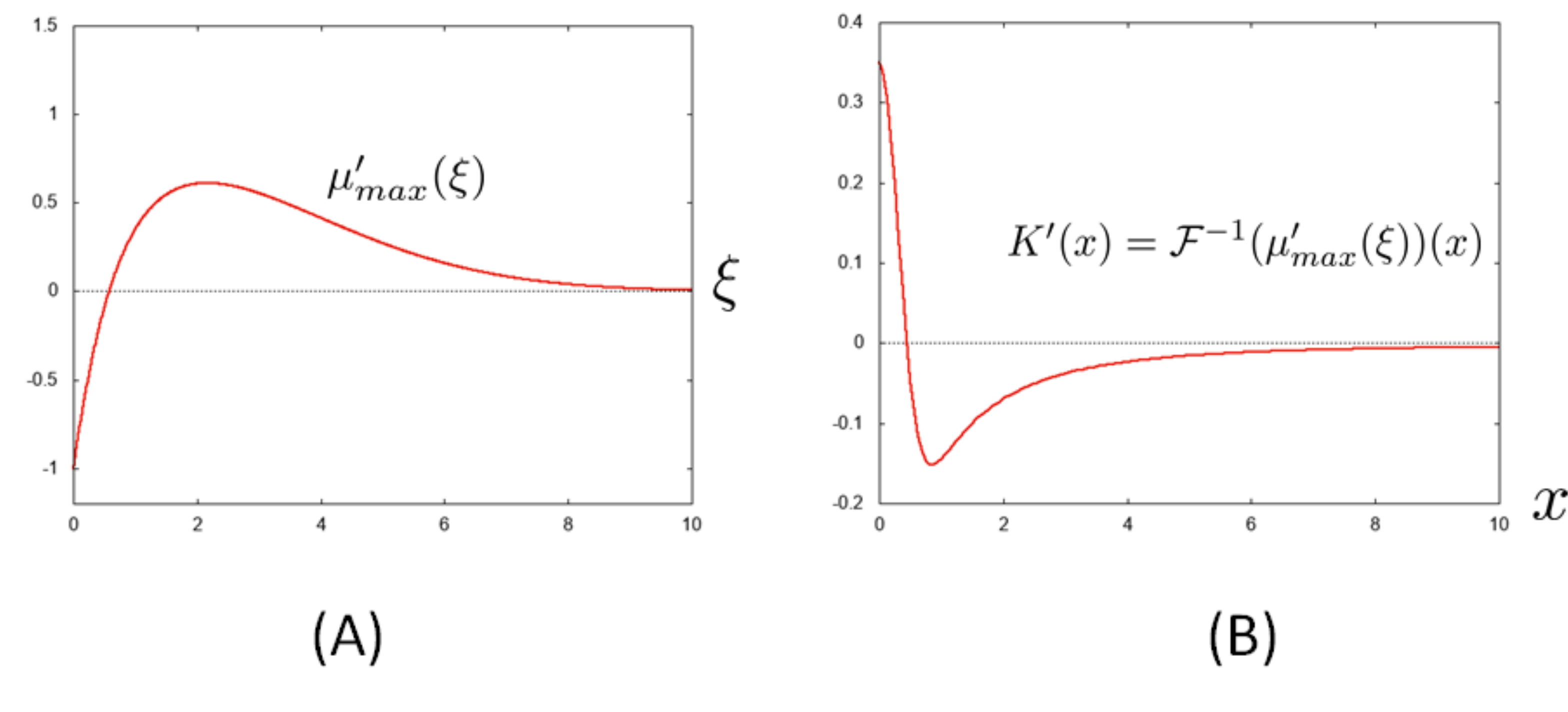}
 \caption{
 (A) Graph of $\mu'_{max}(\xi ) = \dmax_j\{ \zeta'_j(\xi )\}$
for the eigenvalues $\zeta'_j(\xi )$ of $B'_{\vep}(\xi )$ .
 (B) The derived kernel $K'(x) = ({\cal F}^{-1}\mu'_{max})(x)$,
 which shows the Mexican hat profile.
 Parameters are same as those of Fig.\ref{fig7}, Fig.\ref{fig6}
 and Fig.\ref{fig6-1}: 
  $c_1 = c_2 = 1$, $c_3 = 4$, $c_4 = 3$, $d_1 = 0.05$,
 $d_2 = 3$ and  
 $\vep = 0.05$. 
 Only the right half pictures are drawn both in (A) and (B).
 }
 \label{fig6-3}
\end{figure}

When $\zeta_j(\xi )$ is complex satisfying conditions in Section 4.2, we use 
$\wt{\mu}^*_j(\xi ) := M^*(\wt{\omega}_j(\xi ))$ in \rf{e10-41}
as the approximation of $\wt{\mu}_j(\xi )$. Then $\wt{p}(\xi )$,
$\wt{\mu}_j(\xi )$ and $\wt{q}(\xi )$ in \rf{e10-45} are
all calculated as follows:\\
Since 
\begin{eqnarray*}
\wt{\mu}^*_j(\xi ) &=& M^*(\wt{\omega}_j(\xi )) \\
&=&
e^{\delta\omega_j(\xi)}\{
(\omega_j(\xi)\cos\delta b(\xi) - b(\xi)\sin\delta b(\xi) ) \pm
i(b(\xi)\cos\delta b(\xi) + \omega_j(\xi)\sin\delta b(\xi))\}
\end{eqnarray*}
($j=1,2$) as in \rf{e10-41}, 
we see from \rf{e10-21}
\[
\wt{\mu}^*_j(\xi ) = \left\{ \begin{array}{l}
M^*(\zeta_j(\xi ))\; ( \xi < \xi_c), \\
e^{\delta a(\xi)}\{
(a(\xi)\cos\delta b(\xi) - b(\xi)\sin\delta b(\xi) ) \pm
i(b(\xi)\cos\delta b(\xi) + a(\xi)\sin\delta b(\xi))\} \; ( \xi > \xi_c)
\end{array} \right.
\]
and hence we put  for $\xi > \xi_c$
\[
\wt{a}^*(\xi ) := e^{\delta a(\xi)}\{
a(\xi)\cos\delta b(\xi) - b(\xi)\sin\delta b(\xi) \},
\]
\[
\wt{b}^*(\xi ) := e^{\delta a(\xi)}\{
b(\xi)\cos\delta b(\xi) + a(\xi)\sin\delta b(\xi)\},
\]
which leads
\begin{equation} \label{e10-48-1}
\wt{p}^*(\xi ) := \left\{ \begin{array}{clc} 
\frac{\wt{\mu}^*_2(\xi )-\wt{\mu}^*_1(\xi )}{\zeta_2(\xi )-\zeta_1(\xi )}
&= \frac{M^*(\zeta_2(\xi )) - M^*(\zeta_1(\xi ))}{\zeta_2(\xi )-\zeta_1(\xi )}
& (\xi < \xi_c), \\
\wt{p}^*_c &:= e^{\delta a(\xi_c)}(1+\delta a(\xi_c)) & (\xi = \xi_c), \\
\frac{\wt{b}^*(\xi )}{b(\xi )} 
&=  e^{\delta a(\xi)}\{
(\cos\delta b(\xi) + a(\xi)\frac{\sin\delta b(\xi)}{b(\xi )})\}
& (\xi > \xi_c ) \end{array} \right. 
\end{equation}
and
\begin{equation} \label{e10-48-2}
\wt{\nu}^*_j (\xi ) := \left\{ \begin{array}{ll} 
\wt{\mu}^*_j(\xi )  = M^*(\zeta_j(\xi ))  & (\xi \le \xi_c), \\
\wt{a}^*(\xi )  & (\xi > \xi_c ), \end{array} \right. 
\end{equation}
\begin{equation} \label{e10-48-3}
\wt{q}^*(\xi ) := \left\{ \begin{array}{ll} 
0  & (\xi \le \xi_c), \\
-\wt{b}^*(\xi )b(\xi )  & (\xi > \xi_c ). \end{array} \right. 
\end{equation}
Thus, we get the effective equation
\begin{equation} \label{e10-5-1}
\left\{ \begin{array}{ccl}
\dot{X} &=& {\cal L}X + K^* * X + L^**Y, \\
\dot{Y} &=& {\cal L}Y + M^* * X + N^* * Y,
\end{array} \right.
\end{equation}
where $K^*(x) := {\cal F}^{-1}(\wt{\nu}^*_1 (\xi ) )(x)$, 
$L^*(x) := {\cal F}^{-1}(\wt{p}^*(\xi ))(x)$,
$M^*(x) := {\cal F}^{-1}(\wt{q}^*(\xi ))(x)$,
$N^*(x) := {\cal F}^{-1}(\wt{\nu}^*_2(\xi ))(x)$
and the differential operator
${\cal L} := {\cal F}^{-1}(\lambda_h(\xi ))$.

\ \\
Practical way II):\\
In \rf{e2-1}, we use $B'_{\vep}(\xi ) := e^{-\vep\xi^2}B_h(\xi )$
instead of $B_{\vep}(\xi )$. Then all eigenvalues $\zeta'_j(\xi )$ of
$B'_{\vep}(\xi )$ satisfy $\zeta'_j(\xi ) \rightarrow 0$ as $|\xi | \rightarrow \infty$.
 Hence by adjusting $\vep$ and $\delta$ appropriately 
 such that $Re(\zeta'_j(\xi )) > -1/\delta$, we may take
 $\mu = \zeta'_j(\xi )$ as the solution of \rf{e3-3}.
 Actually, $M_{max}(z) = z$ holds for $z \in {\bf C}$ with
 $Re(z) > -1/\delta$ and $| Im(z) | \le C$ for $C > 0$ when
 $0 < \delta << 1$. Consequently, we can simply take
 $\mu'_{max}(\xi ) := \dmax_j\{ Re(\zeta'_j(\xi )\}$
 instead of $\mu_{max}(\xi )$. In particular,
 $\mu'_{max}(\xi ) = \dmax_j\{ \zeta'_j(\xi )\}$ when
 $\zeta'_j(\xi ) \in {\bf R}$.
 In Fig.\ref{fig6-3}, the kernel
$K'(x) = {\cal F}^{-1}(\mu'_{max}(\xi ))(x)$ is drawn under
the same parameters as Fig.\ref{fig6}. It suggests that 
the essential structure of the original kernel $K(x)$ is retained by this rough
approximation.

 In the case of complex $\zeta'_j(\xi )$ as in Section 4.2, 
 we define $b'(\xi )$, $\omega'_j(\xi )$ and $\wt{\omega}'_j(\xi )$
 for $\zeta'_j(\xi )$
 in the same manners as  $b(\xi )$, $\omega_j(\xi )$ and $\wt{\omega}_j(\xi )$
 for $\zeta_j(\xi)$.
 Then, we can take 
 $\wt{\mu}'_j(\xi ) = \wt{\omega}'_j(\xi )$ in \rf{e10-41},
 that is, $\wt{\mu}'_1(\xi ) = \omega'_1(\xi ) + ib'(\xi )$ and
 $\wt{\mu}'_2(\xi ) = \omega'_2(\xi ) - ib'(\xi )$.
Now we can calculate $\wt{p}'(\xi )$, $\wt{\nu}'_j(\xi )$
 and $\wt{q}'(\xi )$ corresponding to
 $\wt{p}(\xi )$, $\wt{\nu}_j(\xi )$  and $\wt{q}(\xi )$
respectively as follows: \\
Since $\wt{\mu}'_j(\xi ) = \zeta'_j(\xi )$ for $\xi < \xi_c$ and
$\wt{\mu}'_1(\xi ) = \ovl{\wt{\mu}'_2(\xi )} = a'(\xi ) + ib'(\xi )$ for $\xi > \xi_c$,
 we see $\wt{a}'(\xi ) = a'(\xi )$, $\wt{b}'(\xi ) = b'(\xi )$ and therefore we get

\begin{equation} \label{e10-50-1}
\wt{p}'(\xi ) := \left\{ \begin{array}{clc} 
\frac{\wt{\mu}'_2(\xi )-\wt{\mu}'_1(\xi )}{\zeta'_2(\xi )-\zeta'_1(\xi )}
&= 1
& (\xi < \xi_c), \\
\wt{p}'_c &:= 1& (\xi = \xi_c), \\
\frac{\wt{b}'(\xi )}{b'(\xi )} 
&=  1
& (\xi > \xi_c ), \end{array} \right. 
\end{equation}
\begin{equation} \label{e10-50-2}
\wt{\nu}'_j (\xi ) := \left\{ \begin{array}{ll} 
\wt{\mu}'_j(\xi )  = \zeta'_j(\xi )  & (\xi \le \xi_c), \\
\wt{a}'(\xi ) = a'(\xi )  & (\xi > \xi_c ), \end{array} \right. 
\end{equation}
and
\begin{equation} \label{e10-50-3}
\wt{q}'(\xi ) := \left\{ \begin{array}{ll} 
0  & (\xi \le \xi_c), \\
-\wt{b}'(\xi )b'(\xi )  = -\{b'(\xi )\}^2 & (\xi > \xi_c ). \end{array} \right. 
\end{equation}
Thus, the effective equation is obtained as
\begin{equation} \label{e10-5-2}
\left\{ \begin{array}{ccl}
\dot{X} &=& {\cal L}X + K' * X + Y, \\
\dot{Y} &=&  {\cal L}Y + M' * X + N' * Y,
\end{array} \right.
\end{equation}
where $K'(x) := {\cal F}^{-1}(\wt{\nu}'_1 (\xi ) )(x)$, 
$M'(x) := {\cal F}^{-1}(\wt{q}'(\xi ))(x)$,
$N'(x) := {\cal F}^{-1}(\wt{\nu}'_2(\xi ))(x)$
and the differential operator
${\cal L} := {\cal F}^{-1}(\lambda_h(\xi ))$.

\SECT{Applications}

In previous sections, we only dealt with 1D problems just for simplicity
while 2D problems should be considered as more realistic situations.
Hence, some of the following applications, we consider 1D and 2D problems.
In 2D case, $\pa_x^2$ is replaced with the Laplacian $\Delta
= \pa_x^2+\pa_y^2$ and
the Heat kernel $H_{\vep}(x)$ on 1D is done with
$H_{\vep}({\bf x}) = \ddfrac{1}{4\pi\vep}
e^{-\frac{r^2}{4\vep}}$ for ${\bf x} = (x,y) \in {\bf R}^2$,
$r = \sqrt{x^2+y^2}$ and so on.

\begin{figure}[h]  
 \centering
 \includegraphics[width=5cm, bb=0 0 676 480]{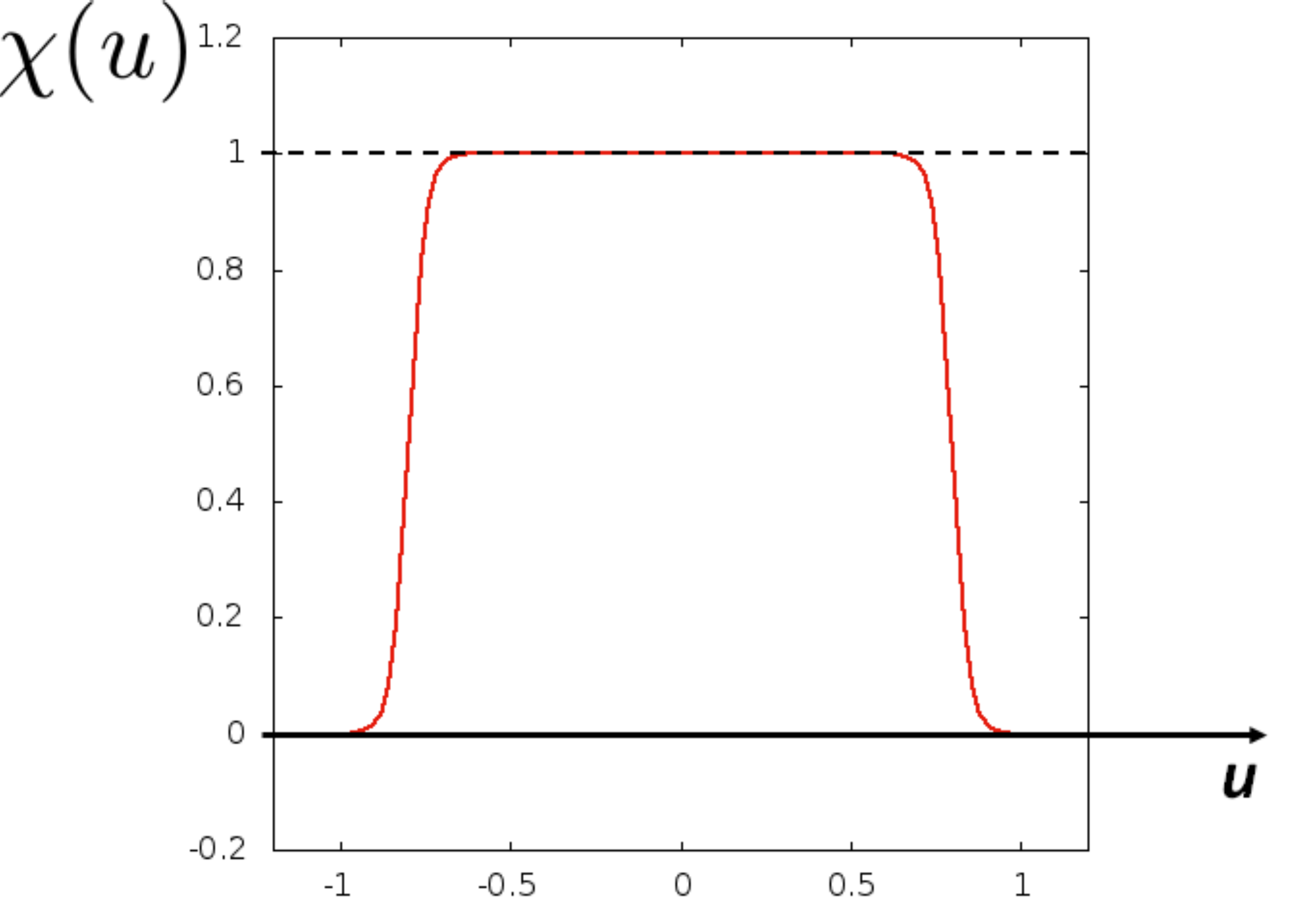}
 \caption{
 An example of the functional form of $\chi (u)$ in
the case of $u^* = 1$.
  }
 \label{fig8-1}
\end{figure}

The numerical simulations on spatial patterns are done by
using the equation
\begin{equation}  \label{e11-1}
u_t = {\cal L}u + \chi (u)\cdot (K_j*u ),
\end{equation}
where $K_j = K_1(x)$ (1D kernel ) or $K_2(x,y)$ (2D kernel)
and ${\cal L}$, the derived kernels and differential operator
in the manner of this paper, and $\chi (u)$ is a cut-off function
satisfying $0 \le \chi (u) \le 1$ and $\chi (u) = 0$ for $|u| \ge u^*$ as in 
Figure \ref{fig8-1}.
%
\subsection{Three node reaction-diffusion network}
\begin{figure}[h]  
 \centering
   \includegraphics[width= 5cm, bb=0 0 484 454 ]{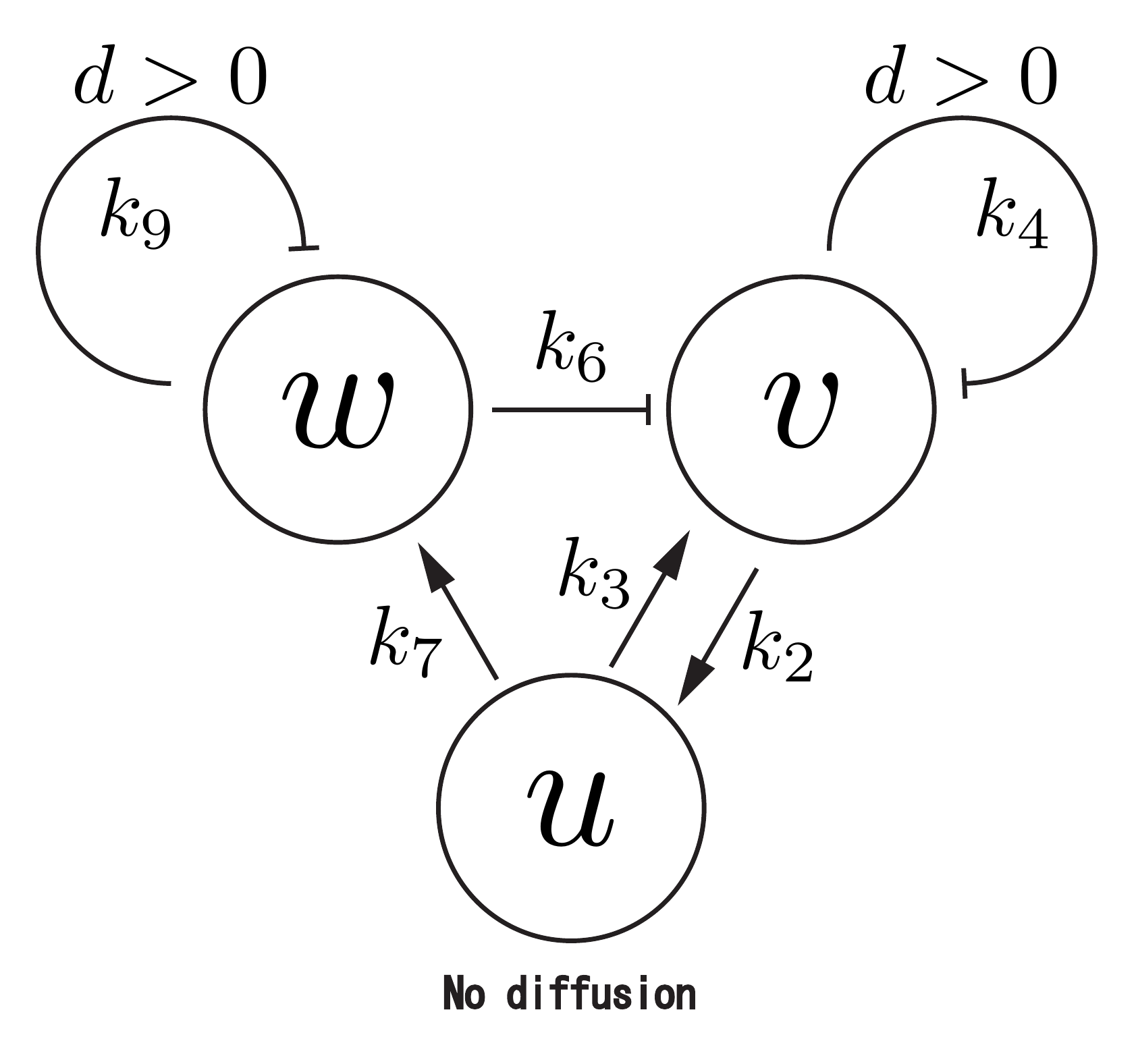}
 \caption{
 Three node network. The picture is cited from \cite{MDSM}
 with the same notations as in Appendix 3-figure 4 in \cite{MDSM}.
  }
 \label{fig8}
\end{figure}

%
In this subsection, we consider the three node network as in Figure \ref{fig8}
dealt with in \cite{MDSM} as Type III network. The typical property of this
example is that $v$ and $w$ diffuse with the same diffusion coefficients $d > 0$
and $u$ does not move. In order to derive an effective kernel of this network,
we use the notations in Section 2. Let $U = {^t}(u,v,w)$. In this network,
the diffusion matrix $D$ is given by $D = diag\{ 0, d, d\}$. 
The matrix $A$ describing
the network is $A = \left( \begin{array}{ccc}
0 & k_2 & 0 \\
k_3 & -k_4 & -k_6 \\
k_7 & 0 & -k_9
\end{array} \right)$. Then the corresponding equation to \rf{e1-1}
leads $\wh{U}_t = B(\xi^2 )\wh{U}$ in 1D and 
$\wh{U}_t = B(R^2 )\wh{U}$ in 2D
with $R = \sqrt{\xi^2+\eta^2}$,
where
\begin{eqnarray*}
B(s ) &=& -sD + A \\
&=& 
-s diag\{ 0,d, d\} 
+ A \\
&=&  
\left(
\begin{array}{ccc}
0 & k_2 & 0 \\
k_3 & -ds - k_4 & -k_6 \\
k_7 & 0 & -ds - k_9
\end{array} \right) 
\end{eqnarray*}
for $s = \xi^2$ or $s = R^2$ depending on the spatial dimension.
Letting $\lambda_j(s)$ ($j=1,2,3$) be eigenvalues of $B(s)$,
we see $\lambda_{max}(s) := \dmax_j\{ \lambda_j(s)\}
\rightarrow c$ as $s \rightarrow \infty$ for $c \in {\bf R}$.
Hence by virtue of the way in Section 4.1, 
we take $\lambda_h(s) = 0$ and $B_h(s) = B(s)$.
Since $B_h(s) = -sD + A$, we should deal with
$B_{\vep}(s) = -sD + e^{-\vep s}A$
for $0 < \vep << 1$ by the manner stated in Section 4.1.
But here we use a simpler way stated in Practical way II)
and put  $B'_{\vep}(s) = e^{-\vep s}B_h(s)$.
Then 
$\mu'_{max}(s)$ is given by $\mu'_{max}(s) 
= \dmax_j\{ \zeta'_j(s) \in \sigma (B'_{\vep}(s))\}$.

\begin{figure}[h]  
 \centering
 \includegraphics[width=12cm,  bb=0 0 913 310]{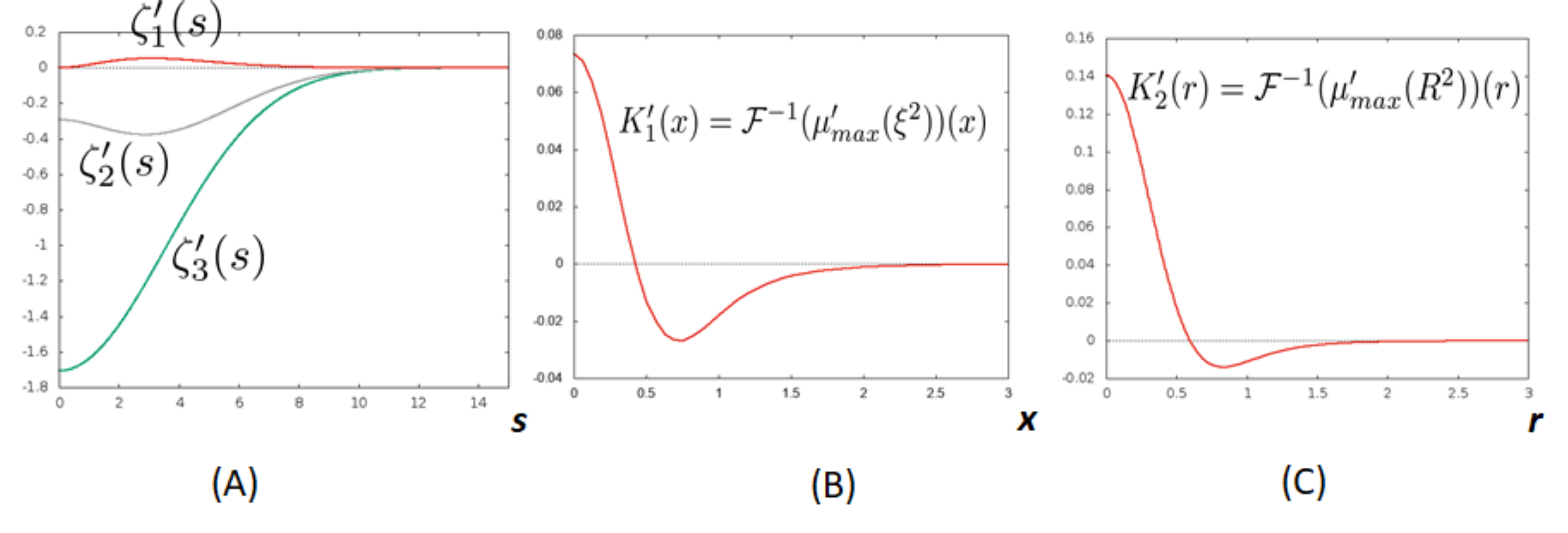}
 \caption{
 (A) Eigenvalues of $B'_{\vep}(s)$, (B) The kernel $K'_1(x) 
 = {\cal F}^{-1}(\mu'_{max}(\xi^2 ))$
 on 1D, (C) The kernel $K'_2(r) = {\cal F}^{-1}(\mu'_{max}(R^2 ))$
 on 2D. 
 Parameter values are 
 $k_2 = 0.5$, $k_3 = k_4 = k_6 = k_7 = k_9 = 1$, $d = 0.02$ and $\vep = 0.05$.
  All values except $\vep$ 
  are same as in Appendix 3-figure 4 in \cite{MDSM}.
  }
 \label{fig9}
\end{figure}

In Fig.\ref{fig9}, the first, second and third maximal eigenvalues of 
$B'_{\vep}(s )$ are drawn and also the 1D and 2D kernels 
$K'_1(x) = {\cal F}^{-1}(\mu'_{max}(\xi^2 ))(x)$
and $K'_2(r) = {\cal F}^{-1}(\mu'_{max}(R^2))(r)$,
which shows that the effective kernel $K'_j(x)$ has a Mexican hat profile and
we can understand the LALI effect essentially occurs even if the diffusion constants
are same.
%

\subsection{Two node reaction-diffusion network with long range interaction}

\begin{figure}
 \centering
 \includegraphics[width=8cm, bb=0 0 811 311]{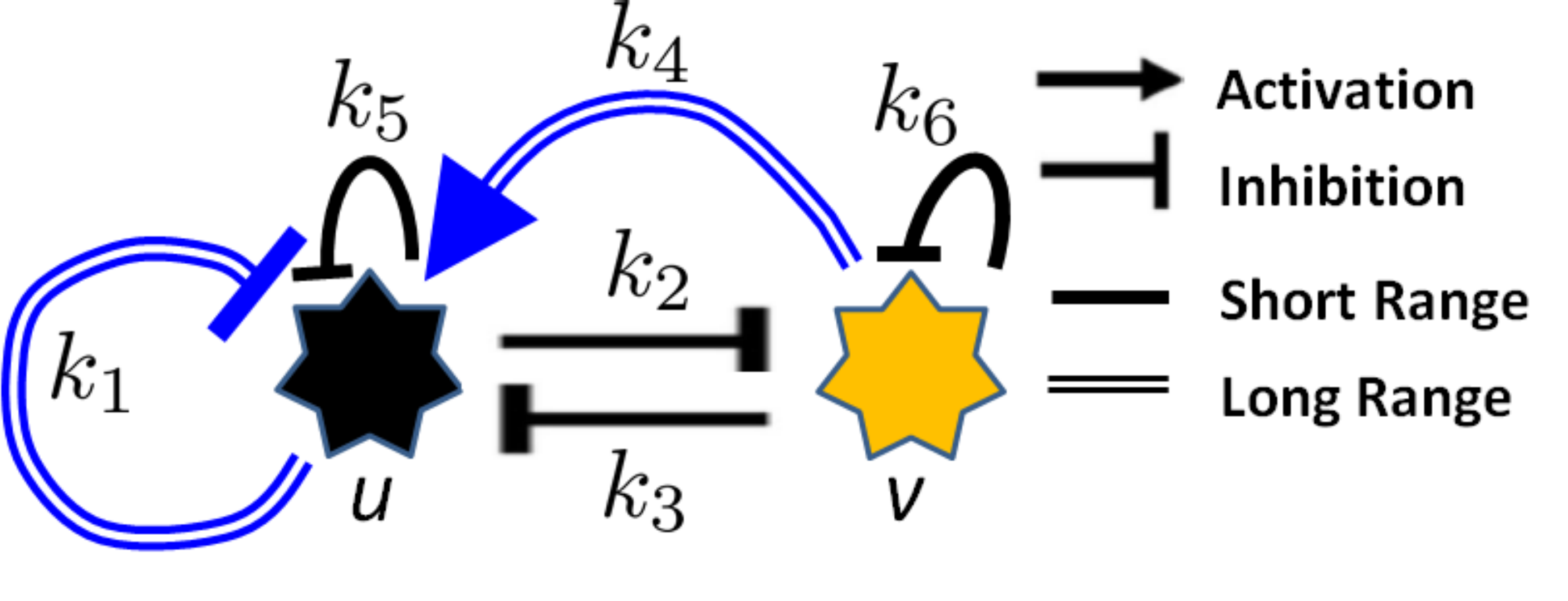}
 \caption{Deduced interaction network between pigment cells of Zebra-fish,
 which is referred from \cite{NTKK,WK}. The coefficient $k_1$ 
is possibly equal to zero (\cite{WK}).
  }
 \label{fig10}
\end{figure}

We apply the method to a problem including long range interaction in the
network, which was stated in \cite{NTKK}, \cite{WK} for interactions 
between zebra-fish pigment cells. The schematic picture of the network
is like Fig.\ref{fig10}.
The long range interaction means a cell interact directly to cells
at a constant distant location by e.g. projections as in Figure \ref{fig4}
and we may represent the effect
by the kernel $L_l(x):= \frac{1}{2}\{
\delta (x+l) + \delta (x-l)\}$ for $x \in {\bf R}$
and $L_l(x,y) 
:= \frac{1}{2\pi l}\delta (|{\bf x}| - l)$ 
for ${\bf x}  = (x,y) \in {\bf R}^2$, where
$l$ is a positive constant.
Let $u$ and $v$ be the differences from some rest states
of densities of the two types of pigment cells,
melanophores and xanthophores, respectively (\cite{NTKK}). 
Then the matrix $A$ describing the network according to Fig.\ref{fig10}
is
\[
AU = \left( \begin{array}{cc}
-k_1L_l* - k_5 & -k_3 + k_4L_l* \\
-k_2 & -k_6
\end{array} \right) U
= \left( \begin{array}{c}
-k_1L_l*u -k_5u - k_3v + k_4L_l*v \\
-k_2u - k_6v
\end{array} \right) 
\]
for $U = {^t}(u,v)$.
In \cite{NTKK}, it was simulated with the same small diffusion constants
of $u$ and $v$. Hence we also assume $u$ and $v$ spatially diffuse
with the same diffusion constant $d$. 
Then the corresponding equation to \rf{e1} is 
\[
U_t = D\Delta U + AU ,
\]
where  $D = diag\{ d, d\}$.
Taking the Fourier transformation of  the above equation, we have
\[
\hat{U}_t = B_1(\xi )\hat{U}\;\; (1D\; case), \;
\hat{U}_t = B_2(R )\hat{U}\;\; (2D\; case),
\]
where $R = \sqrt{\xi^2+\eta^2}$ and
\[
B_j(s) =
\left( \begin{array}{cc}
-ds^2 - k_1\wh{L_l}(s) - k_5 & -k_3 + k_4\wh{L_l}(s) \\
-k_2 & -ds^2 - k_6
\end{array}\right) 
\]
for $s = \xi$ in 1D or $s = R$ in 2D.
Here, $\wh{L_l}$ is computed for 1D, 
\[
\wh{L_l}(\xi ) = \frac{1}{2}(e^{i\xi l}+e^{-i\xi l})\wh{\delta}(\xi )
= \cos\xi l =: P_1(\xi )
\]
and for 2D,
\begin{eqnarray*}\ddfrac{1}{2\pi l}
\wh{L_l}(R) &=& 
\!\dint_0^{\infty}\!\!\!\!\dint_0^{\pi/2}\!\!\!\!\!\!
 r\delta (r-l)\cos ( rR\sin\theta )d\theta dr  
 = \ddfrac{1}{2\pi}\dint_0^{\pi/2}\!\!\!\!\!\!
 \cos ( lR\sin\theta )d\theta dr  
 =: P_2(R)
 \end{eqnarray*}
by using \rf{e11-1}.
Thus we find
\[
B_j(s ) =
\left( \begin{array}{cc}
-ds^2 - k_1P_j(s) - k_5 & \;-k_3 + k_4P_j(s) \\
-k_2 & -ds^2 - k_6
\end{array}\right) 
\]
satisfying $B_j(-s) = B_j(s)$.

Let $\lambda_{1,j}(s)$ and $\lambda_{2,j}(s)$ be eigenvalues
of $B_j(s)$. Since eigenvalues have the asymptotic profile in the highest
order of $s$ as
$\lambda_{i,j}(s) \rightarrow -ds^2 + O(1)$ as $s \rightarrow \infty$,
we  take $\lambda_h(s) = -ds^2$ and $B_{j,h}(s) = 
\{B_j(s)+ds^2I\}$ according to Section 4.1. Since $B_{j,h}(s) = O(1)$,
we put $B_{j,\vep}(s) = e^{-\vep s^2}B_{j,h}(s)$.
Then $\mu_{max}(s)$ may be
calculated by using the approximation stated in Practical way II) as
$\mu_{max}(s) = \dmax_j\{ \zeta_{i,j}(s) \in \sigma (B_{i,\vep}(s)\}$
and we get the effective equation by
\[
(1D):
w_t = dw_{xx} + K_1*w, \;\;
(2D): w_t = d\Delta w + K_2*w.
\]

\begin{figure}[h] 
 \centering
 \includegraphics[width=11cm, bb=0 0 691 540]{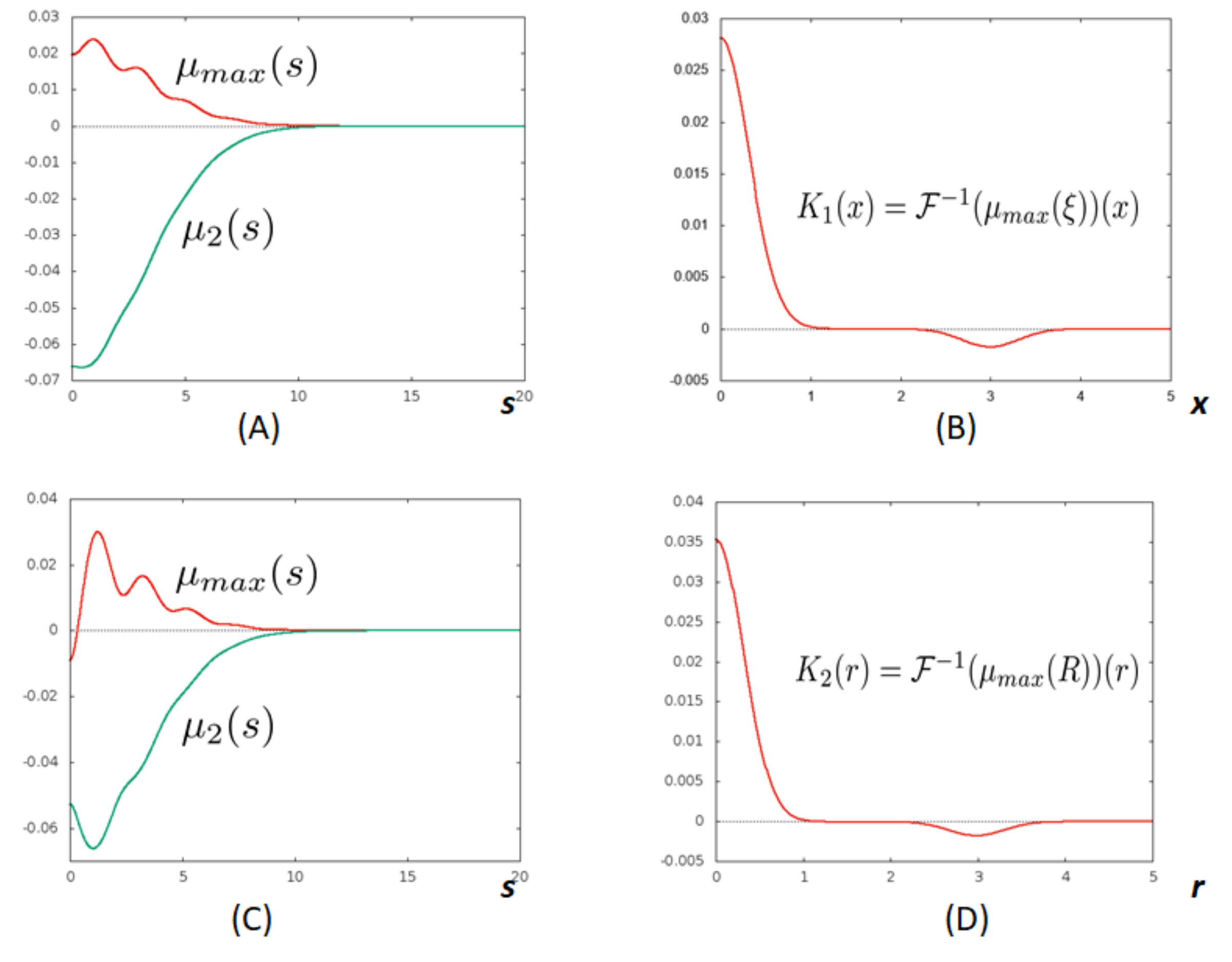}
 \caption{
 (A) $\mu_{max}(s)$ of $B_1(s)$, 
 (B) The kernel $K_1(x) = {\cal F}^{-1}(\mu_{max}(\xi ))(x)$
 by $B_1(s)$,
 (C) $\mu_{max}(s)$ of $B_2(s)$, 
 (D) The kernel $K_2(r) = {\cal F}^{-1}(\mu_{max}(R ))(x)$
 by $B_2(s)$.
 Parameters are
 $d = 0.02$, $k_1 = 0.055*0.016$, $k_2 = 0.05$, $k_3 = 0.04$,
$k_4 = 0.055*0.03$,  $k_5 = 0.02$, $k_6 = 0.025$.
 which are adjusted to parameter values for a wild type of
Fig 4 in \cite{NTKK}. $l$ and $\vep$ are given by $l = 3$
and $\vep = 0.05$.
  }
 \label{fig10-1}
\end{figure}

\begin{figure}[h]  
 \centering
 \includegraphics[width=11cm, bb=0 0 887 328]{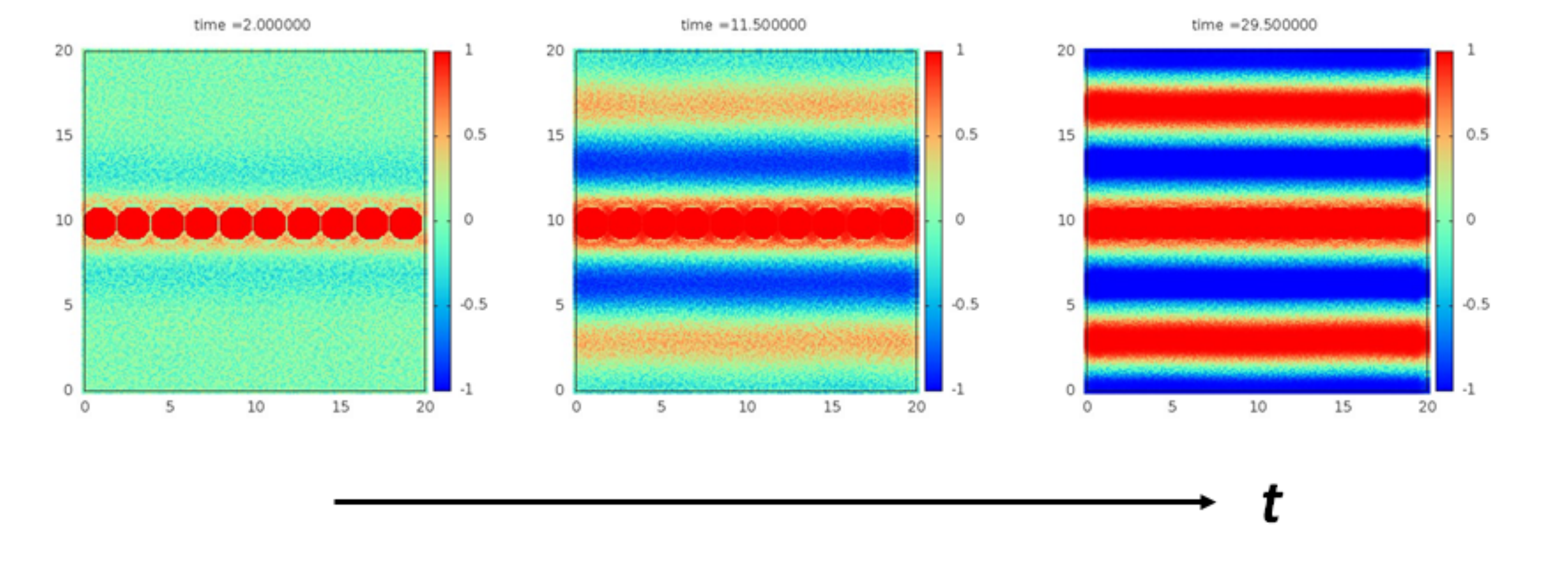}
 \caption{
 Time evolution of spatial patterns in 2D
by the kernel $K_2(r)$ in Fig.\ref{fig10-1} (D) using
\rf{e11-1} with $u^* = 1$.
The initial condition is set similarly to that of the numerical simulations 
reported in \cite{NTKK}.
  }
 \label{fig10-2}
\end{figure}

\begin{figure}
 \centering
 \includegraphics[width=12cm, bb=0 0 933 269]{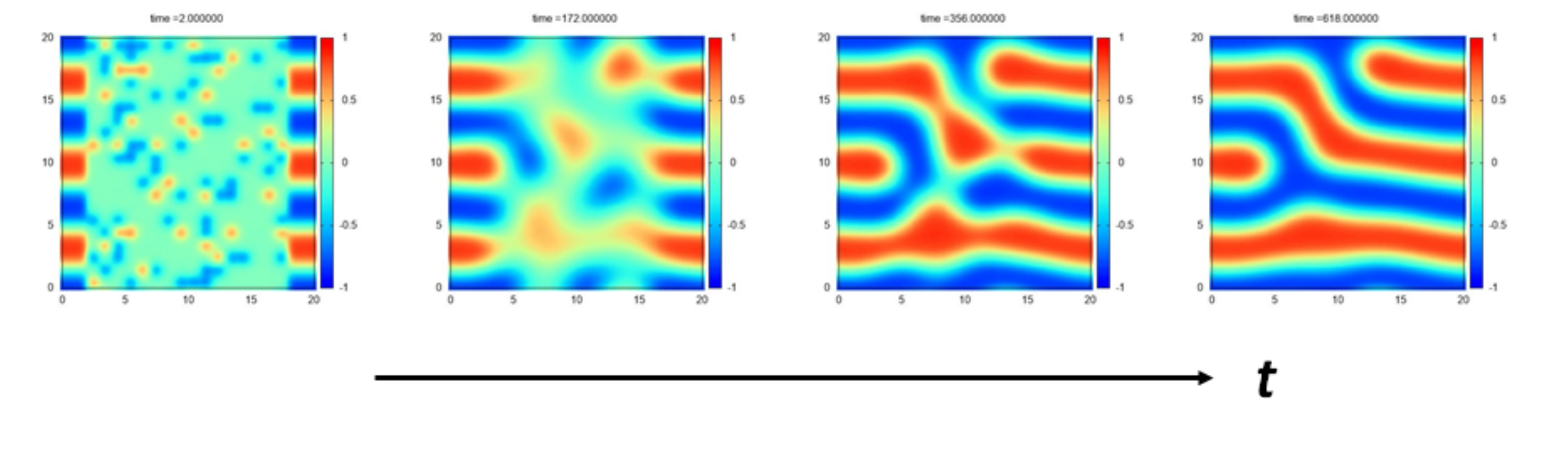}
 \caption{
 Time evolution of spatial patterns in 2D
 with initial data added a random ablated region 
by the kernel $K_2(r)$ in Fig.\ref{fig10-1} (D) using
\rf{e11-1} with $u^* = 1$.
  }
 \label{fig10-3}
\end{figure}

In Fig.\ref{fig10-1},  $\mu_{max}(s)$ of $B_1(s)$ and $B_2(s)$ are respectively drawn 
together with the reduced kernels in 1D and 2D
for the case corresponding to the simulation in \cite{NTKK}.
Fig.\ref{fig10-2} and Fig.\ref{fig10-3} are numerical simulations for 2D patterns
by using the reduced equation $w_t = d\Delta w + \chi (w)\cdot (K_2*w)$ with the
kernel $K_2(r)$ in Fig.\ref{fig10-1},
which shows that the similar
patterns to \cite{NTKK} appear.  
Here, we should note that the system of $u$ and $v$ without
long range interactions (corresponding to the case of $k_1 = k_4 = 0$
in Figure \ref{fig10})  is unstable in the kinetics, that is,
the matrix $\left( \begin{array}{cc} -k_5 & -k_3 \\
-k_2 & -k_6 \end{array}\right)$ has eigenvalues with positive
real parts under the parameter values in Figure \ref{fig10-1}.
In \cite{NTKK}, it was demonstrated that the Turing pattern appears by introducing 
the third component $w$, which plays the role of long range interactions.
In this subsection, the long range interactions are directly introduced by
$L_l$ and the Turing pattern is naturally observed in 2D as in (C) of
Figure \ref{fig10-1}, that is, $\lambda_1(0) < 0$ in (C).
Thus, the method in this paper does not require
any artificial treatment and make the direct understanding of the mechanism
possible through the reduced kernels.

\begin{REM}
In Fig.\ref{fig10-1} (A), we observe $\lambda_1(0) > 0$,
which means that the Turing pattern does not appear
in 1D case for any $l > 0$.
Actually, the functional form of $P_1(\xi ) \sim 2\cos\xi l$ leads
$P_1(0) = 2$ independent of $l$ and the matrix $B_1(0)$ has an positive eigenvalue
under the parameter values in Fig.\ref{fig10-1} while in 2D,
$P_2(0) = 2\pi l$ holds and the maximal eigenvalue $\lambda_1(0)$ of 
$B_2(0)$ can be negative depending
on $l$.
\end{REM}


\begin{figure}
 \centering
 \includegraphics[width=11cm, bb=0 0 702 540]{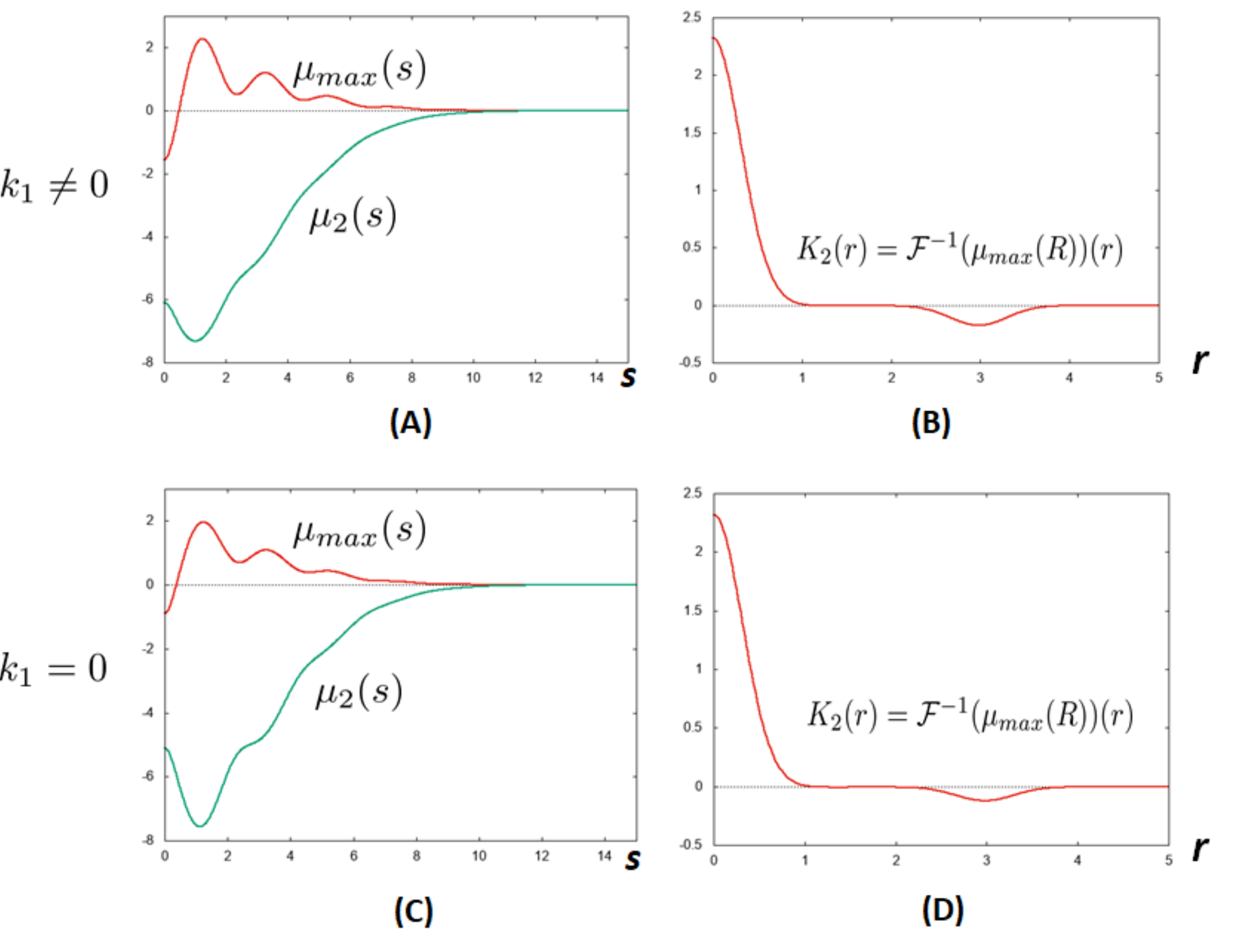}
 \caption{
  Comparison between the cases of  $k_1 = 5.5*0.016 \ne 0$ and $k_1 = 0$.
 Other parameters are fixed by
$l=3.0$, $d=0.2$, $k_2=5.0$, $k_3=4.0$, $k_4=5.5*0.03$,
$k_5=3.0$, $k_6=3.0$. For $k_1 = 5.5*0.05$,
 (A) eigenvalues of $B_2(s)$, 
 (B) the kernel $K_2(r) = {\cal F}^{-1}(\mu_{max}(R))(r)$.  For $k_1 = 0$,
 (C) eigenvalues of $B_2(s)$, 
 (D) the kernel $K_2(r) = {\cal F}^{-1}(\mu_{max}(R ))(r)$.
  }
 \label{fig10-4}
\end{figure}

\begin{figure}
 \centering
 \includegraphics[width=11cm, bb=0 0 852 407]{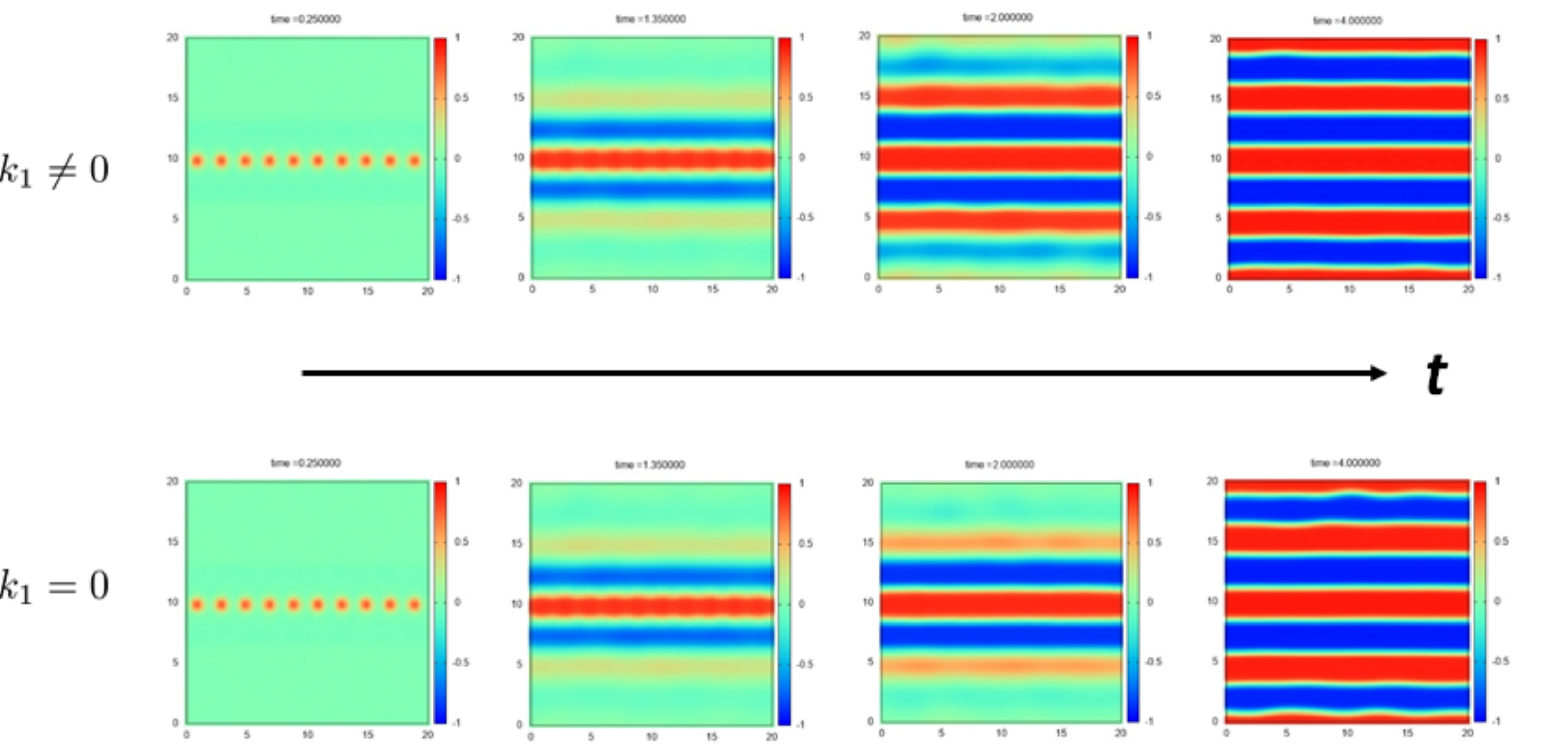}
 \caption{
 Comparison between the cases of  $k_1 = 5.5*0.016 \ne 0$ and $k_1 = 0$.
 Other parameters are fixed under the parameters of Fig.\ref{fig10-4}.
 The upper picture is a time evolutional  2D pattern with the kernel
 $K_2(r)$ of Fig.\ref{fig10-4} (B) and 
 the lower picture is the one with the kernel
 $K_2(r)$ of Fig.\ref{fig10-4} (D).
The initial data is the same one as
 Fig.\ref{fig10-2}.
  }
 \label{fig10-5}
\end{figure}

\begin{figure}
 \centering
 \includegraphics[width=11cm, bb=0 0 837 381]{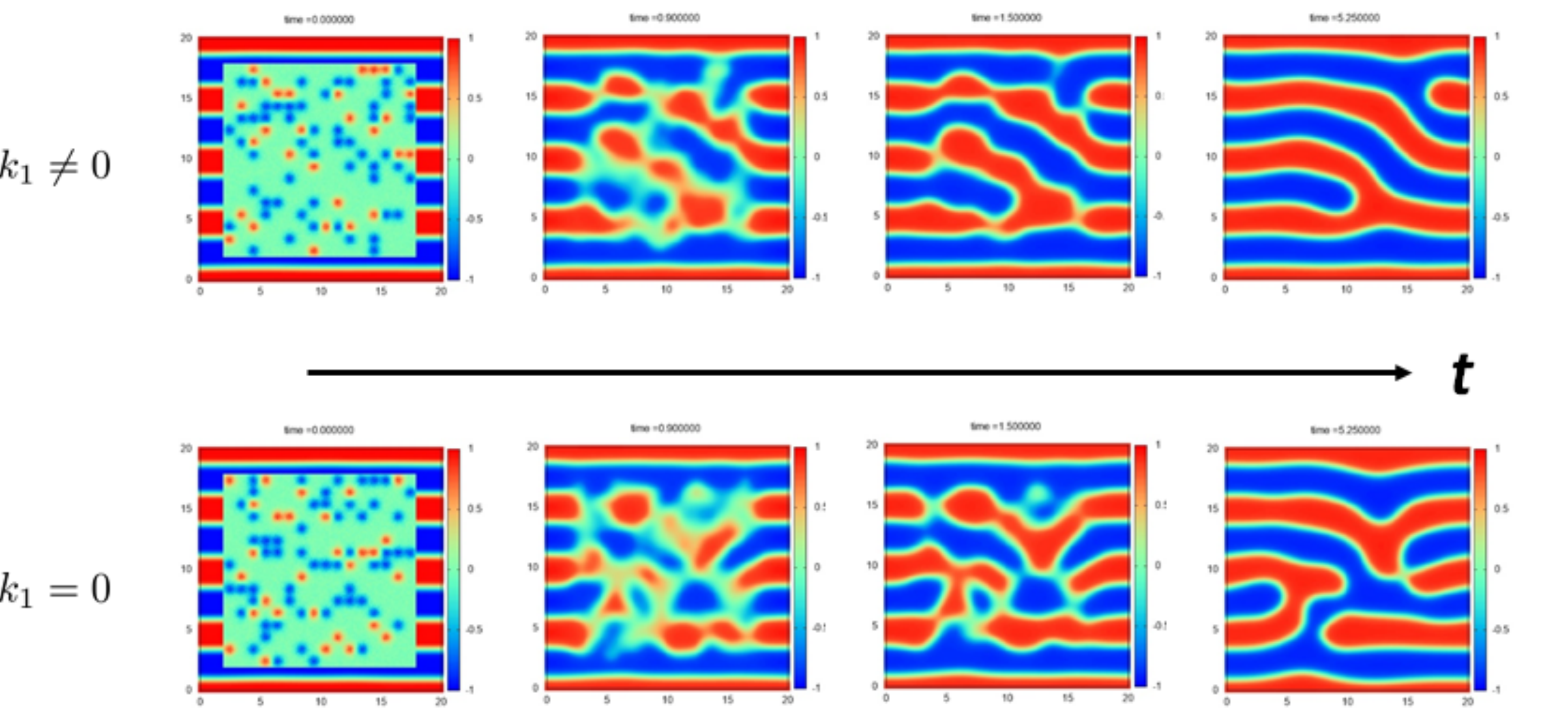}
 \caption{
 Comparison between the cases of  $k_1 = 5.5*0.016 \ne 0$ and $k_1 = 0$.
 Other parameters are fixed under the parameters of Fig.\ref{fig10-4}.
 The upper picture is a time evolutional 2D pattern with the kernel
 $K_2(r)$ of Fig.\ref{fig10-4} (B) and 
 the lower picture is the one with the kernel
 $K_2(r)$ of Fig.\ref{fig10-4} (D).
The initial data is the same one as
 Fig.\ref{fig10-3} added a random ablated region.
  }
 \label{fig10-6}
\end{figure}

In the last of this subsection, we show numerical simulations
comparing two cases of $k_1 \ne 0$ and $k_1 = 0$
because the self-inhibition of $u$ component is not assumed in \cite{WK}.
Fig.\ref{fig10-4} shows the comparison of eigenvalues 
and the reduced kernels between the cases  $k_1 \ne 0$ and
$k_1 = 0$. Fig.\ref{fig10-5} and Fig.\ref{fig10-6}
show the numerical simulations in 2D patterns
under the parameters of Fig.\ref{fig10-4}.
They strongly suggest that
two situations with nonzero and zero $k_1$ values
exhibit almost same properties with respect to 
reduced kernels and 2D patterns. 

\subsection{Network for Proneural waves}
%
%
\begin{figure}
 \centering
 \includegraphics[width=12cm, bb=0 0 886 540]{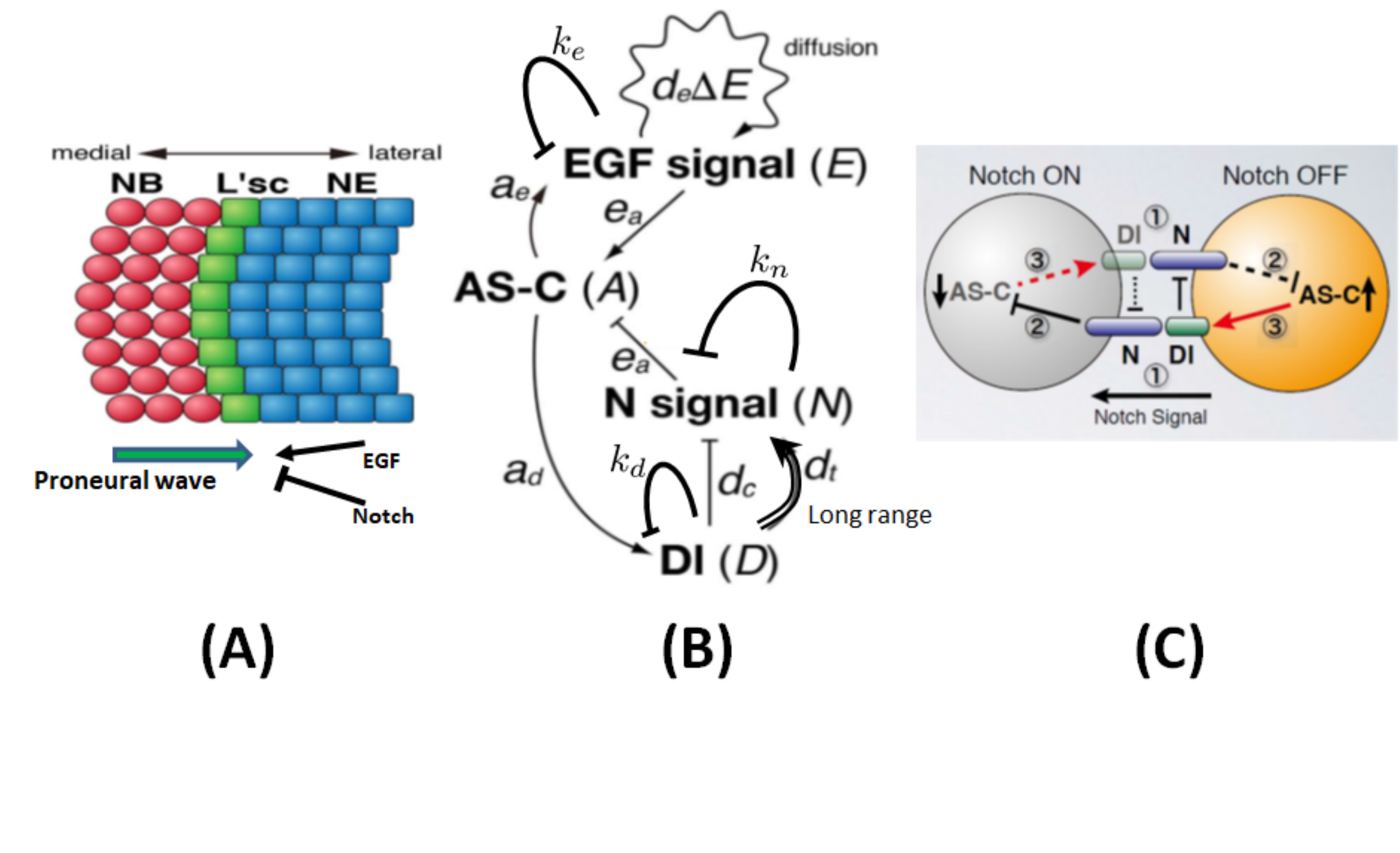}
 \caption{
 The proneural wave progresses unidirectionally during optic lobe development in Drosophila.
 (A) A schematic of the progression of the proneural wave. 
 The proneural wave sweeps from medial to lateral. L'sc is transiently expressed in the 
 differentiating neuroepithelial cells (NEs) and defines the timing 
 of the  differentiation of the NE to the neuroblast (NB). 
 EGF positively regulates the progression of the wave, while Notch negatively regulates 
 wave progression through increasing the expression of Notch target genes.
 (B) Schema showing the gene regulatory network including AS-C, EGF, Notch and Delta. 
 (C) Lateral inhibition mechanism of Delta/Notch signaling.
 }
 \label{fig11}
\end{figure}

\begin{figure}
\begin{center}
\includegraphics[width=11cm, bb=0 0 899 391]{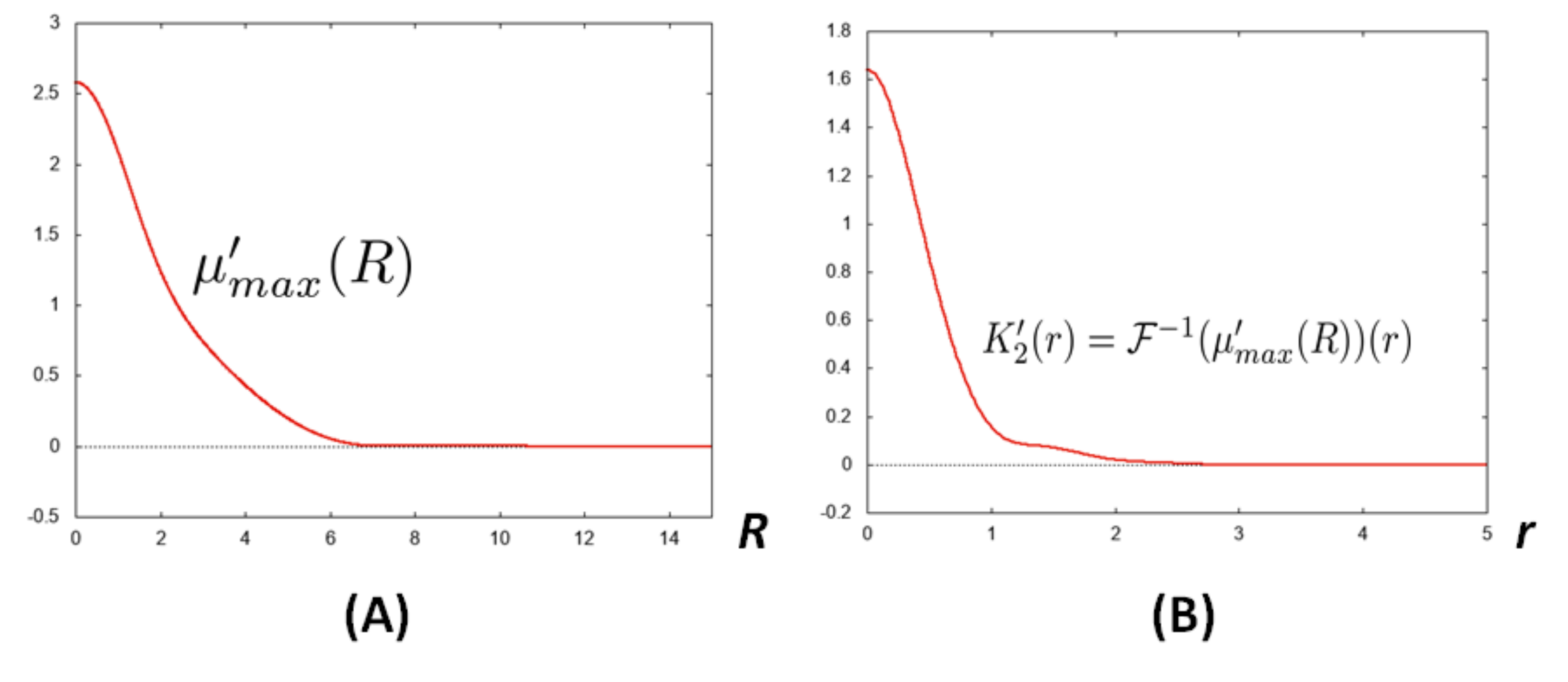}
 \caption{
 (A)~ $\mu'_{max}(R)$ obtained from  $B'_{2,\vep}(s)$. (B)~ 
 The reduced kernel from $\mu'_{max}(R)$. 
 The parameter values are
 $d_e = 1.0$, $k_e = 1.0$, $a_e=1.0$, 
 $k_n=2.0$, $d_t=0.5/(2\pi l)$, $d_c=0.1$,
 $k_d=1.5$, $a_d=1.0$ and $e_a=10.0$.
We also take $l = 1$ and $\eps=0.05$.
  }
 \label{fig11-1}
 \end{center}
\end{figure}

\begin{figure}
 \centering
 \includegraphics[width=12cm, bb=0 0 960 256]{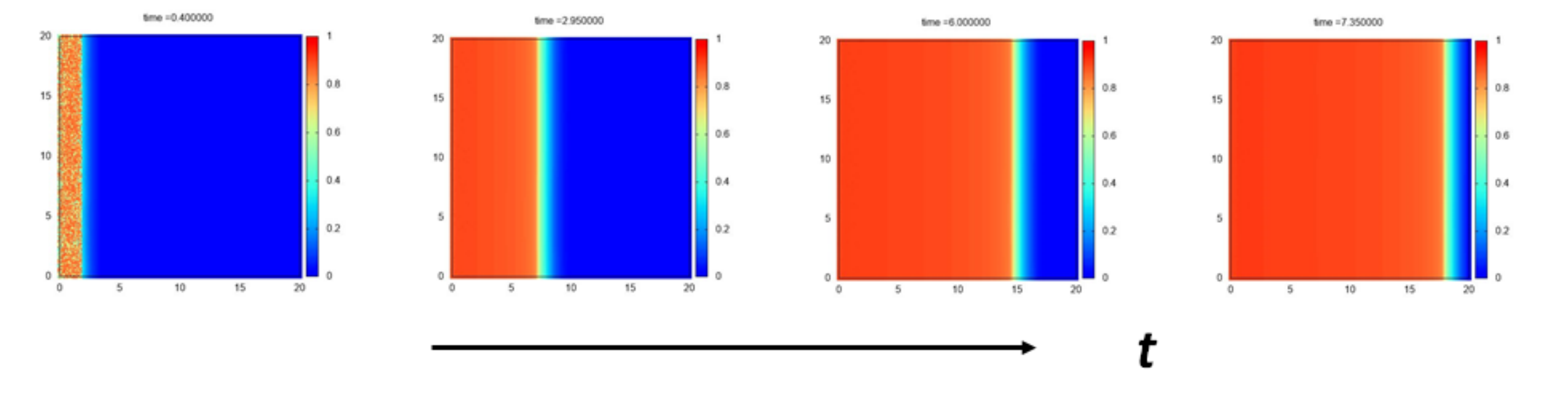}
 \caption{
Time evolution of spatial pattern in 2D simulated by \rf{e5-2}
with the kernel $K'_2(r)$ of Fig.\ref{fig11-1}, which
shows a regulated propagation of proneural wave.
Red color means the higher value of $w$. 
Parameter values are same as those
in Fig.\ref{fig11-1}.
  }
 \label{fig11-2}
\end{figure}

It was reported in \cite{SYMMN} that 
regulated waves of differentiation called the ``proneural wave''
is observed in visual systems on the surface of the brain of the fruit fly.
The phenomena was modeled in a mathematical PDE model as the 
Delta-Notch system including variables $E$, $N$, $D$ and $A_s$ corresponding
to the EGF ligand concentration and EGF signaling ($E$),
Notch signal activity ($N$), Delta expression ($D$)
and the level of the differentiation of AS-C in cells ($A_s$), respectively.
The network of the molecular interactions is as in Fig.\ref{fig11} (B).
The typical properties of the network is that $E$ diffuses and other
materials 
do not and that $D$ inhibits $N$ in the same cell while it activates $N$
in contiguous cells (Fig.\ref{fig11} (C)).
Then by denoting the distance between a cell and the contiguous cells by $l > 0$,
the matrix $A$ describing the network drawn in Fig.\ref{fig11}
is given by
\[
AU = \left(\begin{array}{cccc}
 -k_e & 0 & 0 & a_e \\ 
0 & -k_n & d_tL_l* - d_c & 0 \\
0 & 0 & -k_d & a_d \\ 
e_a & -e_a & 0 & 0
\end{array}\right)U
=
\left(\begin{array}{c}
-k_eE + a_eA_s \\ -k_nN + d_tL_l*D - d_cD \\
-k_dD + a_dA_s \\ e_a(E - N)
\end{array}\right)
\]
for $U = {^t}(E,N,D,A_s)$. Here we used the same functions and notations as those in 
the previous subsections such as $P_j(s)$. The equation is
\[
U_t = D\Delta U + AU ,
\]
where  $D = diag\{ d_e, 0,0,0\}$.
The Fourier transformation leads
$\wh{U}_t = B_1(\xi)\wh{U}$ (1D case) and
$\wh{U}_t = B_2(R)\wh{U}$ (2D case), where $R = \sqrt{\xi^2+\eta^2}$,
$B_j(s) = -s^2D + \wh{A}_j(s)$ and
\[
\wh{A}_j(s) = \left( \begin{array}{cccc}
- k_e & 0 & 0 & a_e \\
0 & -k_n & d_tP_j(s) - d_c & 0 \\
0 & 0 & -k_d & a_d \\
e_a & -e_a & 0 & 0
\end{array} \right) .
\]
Since the maximal eigenvalue of $B_j(s)$ is $O(1)$, we take
$\lambda_h(s) = 0$, $B_{j,h}(s) = B_j(s) = -s^2D + \wh{A}_j(s)$ 
with $\wh{A}_j(s) = O(1)$ and therefore
we should put $B_{j,\vep}(s) = -s^2D + e^{-\vep s^2}\wh{A}_j(s)$
according to the manner in Section 4.
But, here we use the simplest case of 
$B'_{j,\vep}(s) = e^{-\vep s^2}B_{j,h}(s)$ according to 
Practical way II) and only
consider 2 dimensional problems (the case of $s = R = \sqrt{\xi^2+\eta^2}$). 
The effective equations are as follows:

When all eigenvalues $\zeta'_j(R)$ of $B'_{j,\vep}(R)$
are real, the effective equation is
$w_t = K'_2*w$, where $K'_2(r) = {\cal F}^{-1}(\mu'_{max}(R))(r)$
and $\mu'_{max}(R) = \dmax_j\{ \zeta'_j(R)\}$ for $r = \sqrt{x^2+y^2}$.
In numerical simulations, the equation
\begin{equation} \label{e5-2}
w_t = \chi (w)\max\{ \cdot (K'_2*w), 0\}
\end{equation}
is treated
because of the irreversibility of the differentiation of cells, 
which will require the monotonicity in time and
the modification in \rf{e5-2} seems natural.
In Fig.\ref{fig11-1}, $\mu'_{max}(R)$ and the reduced kernel are
drawn in the case that eigenvalues of $B'_{j,\vep}(s)$ are all real.
Fig.\ref{fig11-2} draws the numerical simulation of \rf{e5-2} by
using the kernel $K'_2(r)$ in Fig.\ref{fig11-1}. It shows
a stable traveling planar pattern appears corresponding to the 
proneural wave.

\begin{figure}
\begin{center}
\includegraphics[width=12cm, bb=0 0 808 540]{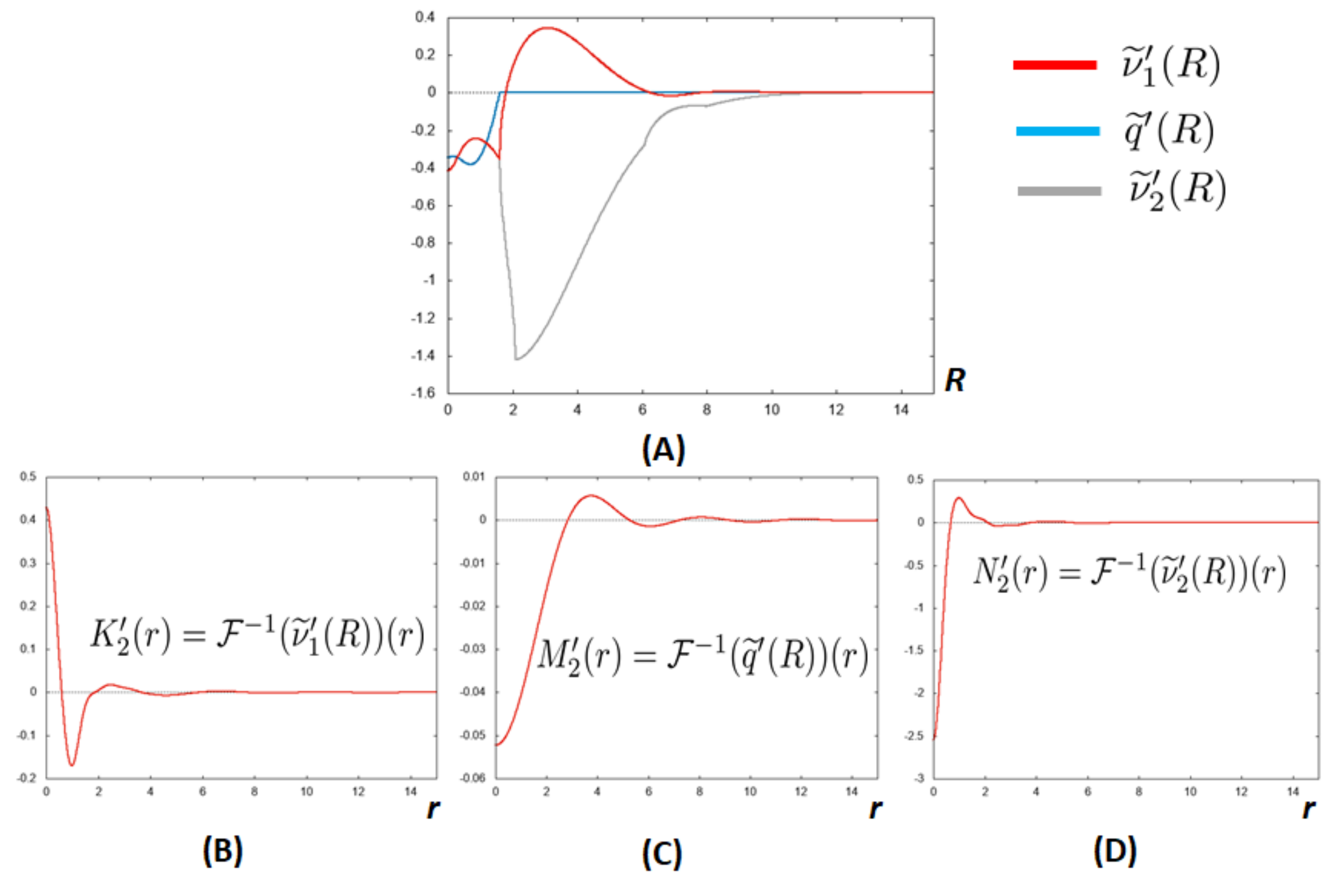}
 \caption{
 (A)~ $\wt{\nu}'_1(R)$, $\wt{\nu}'_2(R)$ and $\wt{q}'(R)$ obtained from 
 $B'_{2,\vep}(R)$. (B), (C), (D) are
 the kernels obtained from them.
 The parameter values are
 $d_e = 1.0$, $k_e = 1.0$, $a_e=0.1$, 
 $k_n=2.0$, $d_t=0.5/(2\pi l)$, $d_c=0.1$,
 $k_d=1.5$, $a_d=1.0$ and $e_a=10.0$.
We also take $l = 1.0$ and $\eps=0.05$. Only $a_e$ is
changed smaller than the one in Fig.\ref{fig11-1} and Fig.\ref{fig11-2}. 
  }
 \label{fig12}
 \end{center}
\end{figure}

\begin{figure}
 \centering
\includegraphics[width=12cm, bb=0 0 908 259]{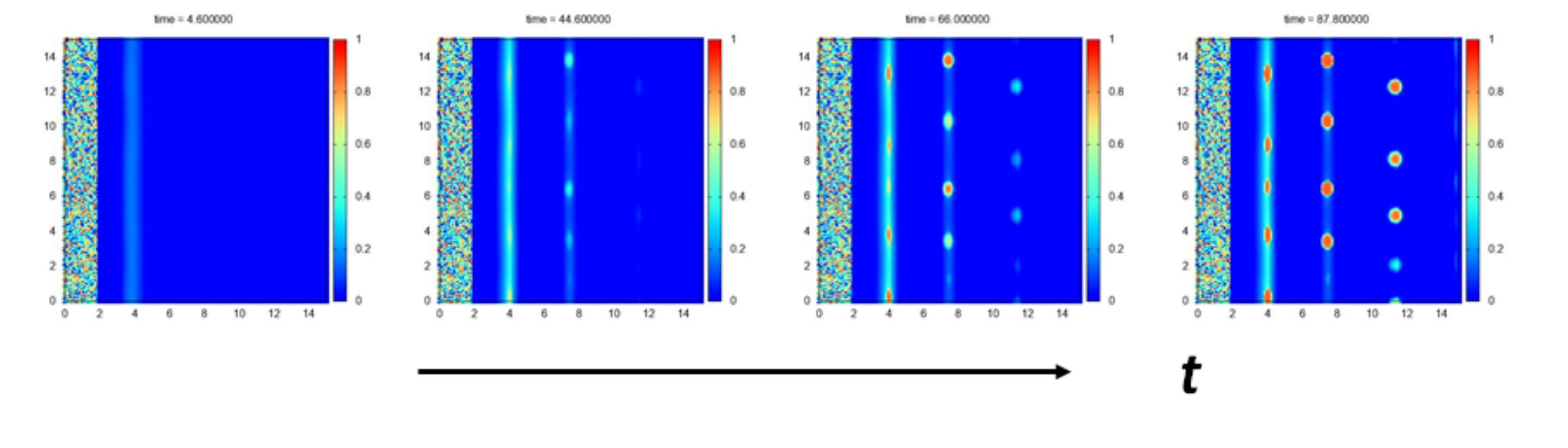}
 \caption{
Time evolution of spatial pattern in 2D simulated by \rf{e5-1}
with kernels $K'_2$, $M'_2$ and $N'_2$ computed in (B), (C)
and (D) of Fig.\ref{fig12}, which shows salt and pepper patterns.
Only $X$ component is drawn. Parameter values are same as those
of Fig.\ref{fig12}. 
  }
 \label{fig13}
\end{figure}

When the eigenvalue with maximal real part can be complex, 
we get \rf{e10-5-2} with ${\cal L} = {\cal F}^{-1}(\lambda_h(R)) = 0$.
that is
\[
\left\{ \begin{array}{ccl}
\dot{X} &=& K'_2 * X + Y, \\
\dot{Y} &=&  M'_2 * X + N'_2 * Y,
\end{array} \right.
\]
where $K'_2(r) = {\cal F}^{-1}(\wt{\nu}'_1(R))(r)$, 
$M'_2(r) = {\cal F}^{-1}(\wt{q}'(R))(r)$ and
$N'_2(r) = {\cal F}^{-1}(\wt{\nu}'_2(R))(r)$.
Numerical simulations are done 
by the following equation:
\begin{equation} \label{e5-1}
\left\{ \begin{array}{ccl}
\dot{X} &=& \chi (X)\max\{ (K'_2 * X + Y), 0\}, \\
\dot{Y} &=& \chi (Y) (M'_2 * X + N'_2 * Y)
\end{array} \right.
\end{equation}
because of the irreversibility of the differentiation of cells
as stated for \rf{e5-2}.
In Fig.\ref{fig12}, the reduced kernels $K'_2(r)$, $M'_2(r)$
and $N'_2(r)$ are shown in the case when $B'_{2\vep}(R)$
can have complex eigenvalues. The parameter values are same
as those of Fig.\ref{fig11-1} except the coefficient $a_e$.
In practice, $a_e = 1$ in Fig.\ref{fig11-1} and $a_e = 0.1$
in Fig.\ref{fig12}. As in Fig.\ref{fig11} and also in the matrix $A$,
$a_e$ denotes the activation rate of $A_s$ (AS-C) for $E$ (EGF).
Lower $a_e$ is expected to enhance the lateral inhibition by Delta/Notch signaling
shown in Fig.\ref{fig11} (C) and as the consequence, pepper and salt
patterns are caused as indicated in \cite{SYMMN}.
Fig.\ref{fig13} demonstrates the occurrence of salt and pepper patterns.

\begin{REM}
In this section, we  adopted the heat kernel $H_{\vep}(x)$
as the approximation of the Dirac $\delta$-function.
If zero state is unstable as the grand state in a system,
the approximation by the heat kernel may cause some trouble
because the heat kernel $H_{\eps}$ activates the instability
everywhere by the property that
$(H_{\eps}*u)(x) > 0$ everywhere for a function $u(x) \ge 0$.
In that case, the approximation by the mollifier 
$J_{\eps}$
stated in Proposition \ref{prop1} seems better
because the mollifier can approximate supports, i.e., the support of
$J_{\eps}*u$ is close to the support of $u$.
Then, we can replace
$B_{j,\vep}(s)$ either by 
$B^*_{j,\vep}(s) = -s^2D + \wh{J_{\vep}}(s)\cdot\wh{A}_j(s)$ or
roughly $B'_{j,\vep}(s) = \wh{J_{\vep}}(s)\cdot B_j(s)
= \wh{J_{\vep}}(s)\cdot\{ -s^2D + \wh{A}_j(s)\}$. 

\end{REM}
\SECT{Discussion}
The efficiency of KT models for the investigation
of spatial patterns in biology was clearly mentioned in \cite{Kon1}.
In this paper, we proposed a method to get the effective kernels
from given network systems. Patterns generated from network systems
are directly  determined by the reduced kernels. It means 
even when given two network systems seem
quite different each other, 
the two different network systems can be identified from the view point
of patterns through the reduced effective kernels. 
The application stated in 5.1 is a typical example of
that point. In practice, the effective kernel of the network of Fig.\ref{fig8}
has a Mexican hat profile as in Fig.\ref{fig9}
with LALI effect similar to the kernel mentioned in Section 3
 and we can regard the generated pattern
as a usual Turing pattern.
Thus we can classify network systems
through effective kernels from the pattern formation point
of view.

Other point of this paper is to
propose a new method for modeling of phenomena.
Two applications given in 5.2 and 5.3 are demonstrations  of it.
In fact, we can observe
the reproduction of the same patterns in Fig.\ref{fig10-2} and
Fig.\ref{fig10-3} 
as the patterns simulated by using a corresponding mathematical model
in \cite{NTKK}
and also the same traveling planar and/or salt-pepper patterns 
in Fig.\ref{fig11-1} and Fig.\ref{fig13}  
as the patterns in \cite{SYMMN}. 
It means that we can propose a new method to make 
mathematical models describing phenomena
in the  types of \rf{e11-1} and \rf{e5-1} by using effective kernels. 

Of course, the model equations in types of 
\rf{e11-1} and \rf{e5-1} have advantages and disadvantages compared with
standard model equations.
One advantage is the systematic and routine derivation of equations for 
given network systems while the construction of model equations in standard
manners requires a deep understanding on the underlying mechanism of patterns
and careful adjustments of nonlinearities.
The possibility of classification of
patterns by kernels is another advantage. The disadvantage is
that the meaning of variables in the equations are not clear
by the big reduction from original systems. Moreover, the equations
of \rf{e11-1} and \rf{e5-1} are reduced from the linearized systems
of given network systems and therefore they can not include
nonlinear effects, which is also disadvantage.
Thus, the model equations in types of 
\rf{e11-1} and \rf{e5-1} and model equations derived in standard manners
should be complementary to each other.

Finally, we mention about how to detect the kernel shape from the observation in
real experiments. 
We apply the Practical way II directly to \rf{e1-1} as the simplest way,
that is, define $B_{\vep}(\xi ) = e^{-\vep\xi^2}B(\xi )$
and $\mu_{max}(\xi ) = \dmax_j\{\zeta_j(\xi ) \in \sigma (B_{\vep}(\xi ) )\}$
by assuming $\zeta_j(\xi ) \in {\bf R}$.
Then the effective equation is 
$u_t = K*u$, where $K(x) = {\cal F}^{-1}(\mu_{max}(\xi ))(x)$. The kernel $K(x)$
is detected by experiments as follows:
Since $\wh{u}(t+\delta ,\xi ) = e^{\delta\wh{K}(\xi )}\wh{u}(t,\xi )
= (1 + \delta\wh{K}(\xi ) + O(\delta^2))\wh{u}(t,\xi )$ holds for any $t > 0$ and
$0 < \delta << 1$, we see $\wh{K}(\xi ) 
= \ddfrac{1}{\delta}\left(\ddfrac{\wh{u}(t+\delta ,\xi )}{\wh{u}(t ,\xi )} - 1\right)
+ O(\delta )$. This implies practically that the ratio of the Fourier transformations of
two profiles, a profile at time $t$ and the one after short time  directly gives the kernel shape
by focusing on an arbitrarily fixed element of the network.

\section*{Acknowledgment}

The authors thank Akiko Nakamasu (Kumamoto University, Japan) for her helpful discussion.




\end{document}